\documentclass[leqno, draft]{article}

\def\average{\mathop{\rm average}}
\def\dist{\mathop{\rm dist}}
\def\im{\mathop{\rm Im}}
\def\re{\mathop{\rm Re}}
\def\span{\mathop{\rm span}}

\newtheorem{theorem}{Theorem}
\newtheorem{lemma}[theorem]{Lemma}
\newtheorem{proposition}[theorem]{Proposition}
\newtheorem{definition}[theorem]{Definition}
\newtheorem{remark}[theorem]{Remark}

\newcommand{\begintheorem}{\addtocounter{equation}{1}\begin{theorem}}
\newcommand{\beginlemma}{\addtocounter{equation}{1}\begin{lemma}}
\newcommand{\beginproposition}{\addtocounter{equation}{1}\begin{proposition}}
\newcommand{\begindefinition}{\addtocounter{equation}{1}\begin{definition}}
\newcommand{\beginremark}{\addtocounter{equation}{1}\begin{remark}}

\begin{document}

\title{Some notes on harmonic and holomorphic functions}

\author{Stephen William Semmes \\
	Rice University \\
	Houston, Texas}

\date{}

\maketitle

\begin{abstract}
These notes are concerned with harmonic and holomorphic functions on
Euclidean spaces, where ``holomorphic'' refers to ordinary complex
analysis in dimension $2$ and generalizations using quaternions and
Clifford algebras in higher dimensions.  Among the principal themes
are weak solutions, the mean-value property, and subharmonicity.
\end{abstract}

\tableofcontents

\section{Euclidean spaces}
\label{section on Euclidean spaces}
\setcounter{equation}{0}

\subsection{Real and complex numbers}
\label{subsection on real and complex numbers}

	Let ${\bf R}$ denote the real numbers, equipped with the usual
arithmetic operations and ordering.  We write ${\bf Z}$ and ${\bf
Z}_+$ for the integers and the positive integers, respectively.
The complex numbers are denoted ${\bf C}$, and of course a complex
number $c$ can be written as $a + i \, b$, where $a$, $b$ are real
numbers and $i^2 = -1$.

	Recall that a real number $b$ is said to be an \emph{upper
bound}\index{upper bound for a subset of the real numbers} of a subset
$A$ of ${\bf R}$ if $a \le b$ for all $a \in A$.  If $A$ is a set of
real numbers and $c$ is a real number, then $c$ is said to be the
\emph{least upper bound}\index{least upper bound of a set of real
numbers} or \emph{supremum}\index{supremum of a set of real numbers}
of $A$ if $c$ is an upper bound for $A$ and if $c \le b$ for every
real number $b$ which is an upper bound for $A$.  It is clear from the
definition that the supremum of $A$ is unique if it exists.

	Similarly, a real number $b'$ is said to be a \emph{lower
bound}\index{lower bound for a set of real numbers} for a subset $A$
of ${\bf R}$ if $b' \le a$ for every element $a$ of $A$.  If $A$ is a
subset of ${\bf R}$ and $c'$ is a real number, then $c'$ is said to be
the \emph{greatest lower bound}\index{greatest lower bound of a set of
real numbers} or \emph{infimum}\index{infimum of a set of real
numbers} if $b' \le c'$ for every real number $b'$ which is a lower
bound for $A$.  Again, it is easy to see from the definition that the
infimum is unique if it exists.

	The completeness property of the real numbers, with respect to
the ordering on the real numbers, states that every nonempty subset of
${\bf R}$ which has an upper bound also has a least upper bound.  This
is equivalent to saying that every nonempty subset of ${\bf R}$ which
has a lower bound also has a greatest lower bound.  The equivalence
between the two statements can be shown using multiplication by $-1$,
or by taking the supremum of the set of lower bounds of a subset of
${\bf R}$ to get a greatest lower bound and the infimum of the set of
upper bounds of a subset of ${\bf R}$ to get a least upper bound, as
appriopriate.

	The supremum and infimum of a subset $A$ of ${\bf R}$ are
denoted
\begin{equation}
	\sup A, \quad \inf A,
\end{equation}
respectively, assuming that they exist.  If $f(x)$ is a real-valued
function on some nonempty set $E$, then we may write
\begin{equation}
	\sup_{x \in E} f(x), \quad \inf_{x \in E} f(x)
\end{equation}
for the supremum and infimum of $f(x)$ on $E$, assuming that they
exist.  More precisely, these are the same as the supremum and infimum
of the set
\begin{equation}
	\{f(x) : x \in E\}
\end{equation}
of values of $f$.

	If $x$ is a real number, then the \emph{absolute
value}\index{absolute value of a real number} of $x$ is denoted $|x|$
and defined to be equal to $x$ when $x \ge 0$ and to $-x$ when $x \le
0$.  One can check that
\begin{equation}
	|x + y| \le |x| + |y|
\end{equation}
for all real numbers $x$, $y$, which is the triangle
inequality\index{triangle inequality for real numbers} for the
absolute value, and that
\begin{equation}
	|x \, y| = |x| \, |y|,
\end{equation}
which is to say that the absolute value of a product is equal to the
product of the individual absolute values.  As a consequence of the
triangle inequality, one can check that
\begin{equation}
	\Bigl| |x| - |y| \Bigr| \le |x - y|
\end{equation}
for all $x, y \in {\bf R}$.

	The \emph{complex conjugate}\index{complex conjugate of a
complex number} of a complex number $c = a + i \, b$, $a, b \in {\bf
R}$, is denoted $\overline{c}$ and defined by
\begin{equation}
	\overline{c} = a - i \, b.
\end{equation}
One can check that
\begin{equation}
	\overline{\alpha + \beta}
		 = \overline{\alpha} + \overline{\beta}
\end{equation}
and
\begin{equation}
	\overline{\alpha \, \beta}
		= \overline{\alpha} \, \overline{\beta}
\end{equation}
for all $\alpha, \beta \in {\bf C}$.  If $c \in {\bf C}$, $c = a + i
\, b$, with $a, b \in {\bf R}$, then $a$ and $b$ are called the
\emph{real}\index{real part of a complex number} and
\emph{imaginary}\index{imaginary part of a complex number} parts of
$c$ and are denoted $\re c$, $\im c$, respectively, and we have that
\begin{equation}
	\re c = \frac{c + \overline{c}}{2}, \quad
		\im c = \frac{c - \overline{c}}{2 i}.
\end{equation}

	If $c = a + i \, b$ is a complex number, $a, b \in {\bf R}$,
then the \emph{absolute value}\index{absolute value of a complex
number} or \emph{modulus}\index{modulus of a complex number} is
denoted $|c|$ and defined by
\begin{equation}
	|c| = \sqrt{a^2 + b^2}.
\end{equation}
Notice that
\begin{equation}
	|\re c|, \ |\im c| \le |c|
\end{equation}
and
\begin{equation}
	|c|^2 = c \, \overline{c}.
\end{equation}
It follows that the modulus of a product is equal to the product of
the moduli, so that
\begin{equation}
	|\alpha \, \beta| = |\alpha| \, |\beta|
\end{equation}
for $\alpha, \beta \in {\bf C}$.

	A basic result states that the triangle inequality holds for
complex numbers,\index{triangle inequality for complex numbers} i.e.,
\begin{equation}
	|\alpha + \beta| \le |\alpha| + |\beta|
\end{equation}
for all complex numbers $\alpha$, $\beta$.  Indeed, 
\begin{eqnarray}
	|\alpha + \beta|^2 
   & = & |\alpha|^2 + 2 \re \alpha \, \overline{\beta} + |\beta|^2 \\
   & \le & |\alpha|^2 + 2 |\alpha| \, |\beta| + |\beta|^2  
					\nonumber \\
   & = & (|\alpha| + |\beta|)^2.	\nonumber
\end{eqnarray}
As in the case of real numbers, it follows from the triangle
inequality that
\begin{equation}
	\Bigl| |\alpha| - |\beta| \Bigr| \le |\alpha - \beta|
\end{equation}
for $\alpha, \beta \in {\bf C}$.

	Let $n$ be a positive integer.  We write ${\bf R}^n$ for the
space of $n$-tuples of real numbers, also known as $n$-dimensional
Euclidean space.  If $x = (x_1, \ldots, x_n)$ and $y = (y_1, \ldots,
y_n)$ are elements of ${\bf R}^n$ and $t$ is a real number, then
$x + y$ and $t \, x$ are defined as elements of ${\bf R}^n$ using
coordinatewise addition and scalar multiplication, so that
\begin{equation}
	x + y = (x_1 + y_1, \ldots, x_n + y_n)
\end{equation}
and
\begin{equation}
	t \, x = (t \, x_1, \ldots, t \, x_n).
\end{equation}

	The \emph{inner product}\index{inner product on ${\bf R}^n$}
of $x, y \in {\bf R}^n$ is denoted $\langle x, y \rangle$ and defined
by
\begin{equation}
	\langle x, y \rangle = \sum_{j = 1}^n x_j \, y_j.
\end{equation}
Notice that
\begin{equation}
	\langle x, y \rangle = \langle y, x \rangle
\end{equation}
for all $x, y \in {\bf R}^n$.  If $x, x', y \in {\bf R}^n$ and $t, t'
\in {\bf R}$, then
\begin{equation}
	\langle t \, x + t' \, x', y \rangle
		= t \, \langle x, y \rangle + t' \, \langle x', y \rangle,
\end{equation}
and of course there is a similar linearity property in $y$ by
symmetry.

	The Euclidean norm $|x|$ of $x \in {\bf R}^n$ is defined by
\begin{equation}
	|x| = \biggl(\sum_{j=1}^n x_j^2 \biggr)^{1/2}.
\end{equation}
This is the same as
\begin{equation}
	|x|^2 = \langle x, x \rangle.
\end{equation}
Note that if we identify a complex number $c$ with the element of
${\bf R}^2$ given by $(\re c, \im c)$, then the modulus of $c$ as a
complex number is the same as the Euclidean norm of the corresponding
element of ${\bf R}^2$.

	An important feature of the inner product on ${\bf R}^n$ is
the \emph{Cauchy--Schwarz inequality},\index{Cauchy--Schwarz
inequality} which states that
\begin{equation}
	|\langle x, y \rangle| \le |x| \, |y|
\end{equation}
for $x, y \in {\bf R}^n$.  To see this, one can observe that for each
real number $t$,
\begin{eqnarray}
	0 & \le & \langle x + t \, y, x + t \, y \rangle	\\
	  & = & |x|^2 + 2 \, t \, \langle x, y \rangle + t^2 \, |y|^2.
							\nonumber
\end{eqnarray}
The Cauchy--Schwarz inequality follows by choosing $t$ appropriately.

	For each $x \in {\bf R}^n$ and each $t \in {\bf R}$ we have
that
\begin{equation}
	|t \, x| = |t| \, |x|,
\end{equation}
as one can easily check from the definitions.  Note that $|t|$ refers
to the absolute value of the real number $t$, while $|t \, x|$, $|x|$
refer to the norms of the vectors $t \, x$, $x$ in ${\bf R}^n$.  We
also again have the triangle inequality
\begin{equation}
	|x + y| \le |x| + |y|
\end{equation}
for $x, y \in {\bf R}^n$, since
\begin{eqnarray}
	|x + y|^2 & = & \langle x + y, x + y \rangle		\\
		& = & |x|^2 + 2 \, \langle x, y \rangle + |y|^2
						\nonumber \\
		& \le & |x|^2 + 2 \, |x| \, |y| + |y|^2
						\nonumber \\
		& = & (|x| + |y|)^2.		\nonumber
\end{eqnarray}

	Two other important features of the Euclidean norm on ${\bf
R}^n$ are the \emph{polarization identity}\index{polarization identity
for the Euclidean norm on ${\bf R}^n$} and the \emph{parallelogram
law}.\index{parallelogram law for the Euclidean norm on ${\bf R}^n$}
The polarization identity expresses the inner product in terms of the
norm, through the formula
\begin{equation}
	\langle x, y \rangle = \frac{1}{2} 
			\biggl(|x + y|^2 - |x|^2 - |y|^2 \biggr)
\end{equation}
for $x, y \in {\bf R}^n$.  The parallelogram law states that
\begin{equation}
	|x + y|^2 + |x - y|^2 = 2 \, |x|^2 + 2 \, |y|^2,
\end{equation}
which relates the lengths of the sides of a parallelogram to the
lengths of the diagonals.

\subsection{Linear algebra in ${\bf R}^n$}
\label{subsection on linear algebra in R^n}

	Let $n$ be a positive integer.  In this subsection we review
matters concerning linear independence, spanning vectors, linear
subspaces, and linear mappings in ${\bf R}^n$.  We begin with a
preliminary fact about systems of linear equations.

	Suppose that $l$ and $m$ are positive integers, and that
for each $j = 1, \ldots, l$ and $k = 1, \ldots, m$ we have a real
number $a_{j,k}$.  This leads to a system of $l$ homogeneous linear
equations in $m$ unknowns, namely the equations
\begin{equation}
	\sum_{k = 1}^m a_{j,k} \, r_k = 0
\end{equation}
for $j = 1, \ldots, l$.  Here $r_1, \ldots, r_m$ are variables which
run through the real numbers.

	A basic fact is that if $m > l$, then this system of equations
has a nontrivial solution.  In other words, in this case there are
real numbers $r_1, \ldots, r_m$ which satisfy the equations above and
which are not all equal to $0$.  This is not difficult to show, by
using the equations to solve for one variable in terms of the other
variables and substituting this into the remaining equations to 
systematically reduce the number of equations.

	A collection of vectors $v_1, \ldots, v_m$ in ${\bf R}^n$ is
said to be \emph{linearly independent}\index{linearly independent
vectors in ${\bf R}^n$} if for any choices of real numbers $r_1,
\ldots, r_m$ we have that
\begin{equation}
	\sum_{j = 1}^m r_j \, v_j = 0
\end{equation}
in ${\bf R}^n$ implies that 
\begin{equation}
	r_1 = \cdots r_m = 0.
\end{equation}
This is equivalent to saying that a vector $x$ in ${\bf R}^n$ can
be expressed in at most one way as a sum of the form
\begin{equation}
	\sum_{j=1}^m s_j \, v_j
\end{equation}
for some real numbers $s_1, \ldots, s_m$.  From the result mentioned
in the previous paragraph it follows that $m$ is necessarily less than
or equal to $n$ for this to happen.

	Two vectors $v$, $w$ in ${\bf R}^n$ are said to be
\emph{orthogonal}\index{orthogonal vectors in ${\bf R}^n$} if
\begin{equation}
	\langle v, w \rangle = 0,
\end{equation}
and in this case we write $v \perp w$.  A collection of vectors $v_1,
\langle v_m$ is said to be orthogonal if $v_j \perp v_k$ for $1 \le j
\, k \le m$, $j \ne k$.  If $v_1, \ldots, v_m$ are orthogonal vectors
in ${\bf R}^n$ such that $|v_j| = 1$ for each $j$, then $v_1, \ldots,
v_m$ is said to be an \emph{orthonormal}\index{orthonormal collection
of vectors in ${\bf R}^n$} collection of vectors in ${\bf R}^n$.

	A collection of nonzero orthogonal vectors $v_1, \ldots, v_m$
in ${\bf R}^n$ is automatically linearly independent.  It is enough
to check this when $v_1, \ldots, v_m$ are orthonormal, since one can
reduce to this case by multiplying the $v_j$'s by nonzero real numbers
so that they have norm equal to $1$.  In this case, if
\begin{equation}
	x = \sum_{j = 1}^m s_j \, v_j
\end{equation}
for some $x \in {\bf R}^n$ and real numbers $s_1, \ldots, s_m$, then
\begin{equation}
	s_j = \langle x, v_j \rangle
\end{equation}
for each $j$, so that the coefficients $s_j$ are determined by $x$.

	If $w_1, \ldots, w_l$ is any collection of vectors in ${\bf
R}^n$, then a \emph{linear combination}\index{linear combination of a
collection of vectors in ${\bf R}^n$} of $w_1, \ldots, w_l$ is any
vector in ${\bf R}^n$ of the form
\begin{equation}
	\sum_{k = 1}^l r_k \, w_k,
\end{equation}
where $r_1, \ldots, r_l$ are real numbers.  Thus a finite collection
of vectors in ${\bf R}^n$ is linearly independent if each element of
${\bf R}^n$ which can be written as a linear combination of vectors in
the collection can be written as such a linear combination in only one
way.  A collection of vectors $v_1, \ldots, v_m$ in ${\bf R}^m$ is
said to be linearly dependent if it is not linearly independent, and
this is equivalent to saying that one of the $v_j$'s can be written as
a linear combination of the other vectors in the collection.

	By a \emph{linear subspace}\index{linear subspace of ${\bf
R}^n$} of ${\bf R}^n$ we mean a subset $L$ of ${\bf R}^n$ such that $0
\in L$ and $x, y \in L$ implies that $r \, x + s \, y \in L$ for all
$r, s \in {\bf R}$.  If $w_1, \ldots, w_l$ is a collection of vectors
in ${\bf R}^n$, then the \emph{span}\index{span of a collection of
vectors in ${\bf R}^n$} of $w_1, \ldots, w_l$ is denoted
\begin{equation}
	\span (w_1, \ldots, w_l)
\end{equation}
and is the subset of ${\bf R}^n$ consisting of all linear combinations
of $w_1, \ldots, w_l$.  It is easy to see that this is always a linear
subspace of ${\bf R}^n$.

	As a basic example, define $e_1, \ldots, e_n \in {\bf R}^n$ by
saying that $e_j$ has $j$th coordinate equal to $1$ and all other
coordinates equal to $0$.  It is easy to see that $e_1, \ldots, e_n$
is an orthonormal collection of vectors, and is linearly independent
in particular.  Also, the span of $e_1, \ldots, e_n$ is equal to ${\bf
R}^n$.

	If $L$ is a linear subspace of ${\bf R}^n$ and $v_1, \ldots,
v_m$ are linearly independent vectors in ${\bf R}^n$ whose span is
equal to $L$, then we say that $v_1, \ldots, v_m$ form a
\emph{basis}\index{basis for a linear subspace of ${\bf R}^n$} for
$L$.  Thus every linear combination of $v_1, \ldots, v_m$ lies in $L$
in this case, and every element of $L$ can be expressed as a linear
combination of the $v_j$'s in exactly one way.  If the $v_j$'s are
also orthogonal, or orthonormal, then we say that they form an
\emph{orthogonal basis},\index{orthogonal basis for a linear subspace
of ${\bf R}^n$} or \emph{orthonormal basis},\index{orthonormal basis
for a linear subspace of ${\bf R}^n$} of $L$.

	For example, the vectors $e_1, \ldots, e_n$ form an
orthonormal basis of ${\bf R}^n$.  This basis is called the
\emph{standard basis}\index{standard basis for ${\bf R}^n$} for ${\bf
R}^n$.  Of course this basis has $n$ elements.

	In general, if $L$ is a linear subspace of ${\bf R}^n$ and
$v_1, \ldots, v_m$ is a basis for $L$, then the
\emph{dimension}\index{dimension of a linear subspace of ${\bf R}^n$}
of $L$ is defined to be $m$.  One can check that the dimension of a
subspace is independent of the choice of basis.  Indeed, if $w_1,
\ldots, w_p$, $v_1, \ldots, v_m$ are vectors in ${\bf R}^n$ such that
each $v_j$ lies in the span of the $w_k$'s and the $v_j$'s are
linearly independent, then $m \le p$.

	Let us note that every linear subspace $L$ of ${\bf R}^n$ has
a basis.  In the case where $L = \{0\}$, we interpret this as using
the ``empty basis'', and the dimension is equal to $0$.  In general
one can get a basis by first choosing a nonzero vector in $L$, and
systematically adding vectors to the collection which are not in the
span of the vectors already selected, until $L$ is equal to the span
of the vectors selected and one gets a basis.

	If $L$ is a linear subspace of ${\bf R}^n$ and $v_1, \ldots,
v_k$ is a collection of linearly independent vectors in $L$, then one
can add vectors to this collection if necessary to get a basis for
$L$.  This follows from the same kind of argument as the one for
showing that there is a basis for each linear subspace of ${\bf R}^n$.
As a consequence, if $L_1$, $L_2$ are linear subspaces of ${\bf R}^n$
such that $L_1 \subseteq L_2$, then the dimension of $L_1$ is less
than or equal to the dimension of $L_2$.

	Now suppose that $L$ is a linear subspace of ${\bf R}^n$ which
is equal to the span of a collection of vectors $w_1, \ldots, w_m$.
In this case the dimension of $L$ is less than or equal to $m$.
Indeed, either $w_1, \ldots, w_m$ are linearly independent, and hence
a basis of $L$, or one can remove vectors from this collection without
changing the span to get a subcollection which is a basis for $L$.

	Suppose that $L_1$, $L_2$ are linear subspaces of ${\bf R}^n$
such that 
\begin{equation}
	L_1 \cap L_2 = \{0\}.
\end{equation}
Given a collection of linearly independent vectors in $L_1$ and
another collection of linearly independent vectors in $L_2$, one can
combine them to get a collection of vectors in ${\bf R}^n$ which is
also linearly independent.  Thus the sum of the dimensions of $L_1$
and $L_2$ is less than or equal to $n$.

	For any two nonempty subsets $A_1$, $A_2$ of ${\bf R}^n$, put
\begin{equation}
	A_1 + A_2 = \{x \in {\bf R}^n : x = a_1 + a_2
				\hbox{ for some } a_1 \in A_1,
						a_2 \in A_2\}.
\end{equation}
If $A_1$, $A_2$ are linear subspaces of ${\bf R}^n$, then so is the
sum $A_1 + A_2$.  If $A_1$, $A_2$ are linear subspaces of ${\bf R}^n$
such that $A_1 + A_2 = {\bf R}^n$, then one can check that the sum of
the dimensions of $A_1$, $A_2$ is greater than or equal to $n$.

	By a \emph{linear mapping}\index{linear mapping on ${\bf
R}^n$} on ${\bf R}^n$ we mean a mapping $T$ from ${\bf R}^n$ to itself
such that
\begin{equation}
	T(x + y) = T(x) + T(y)
\end{equation}
for all $x, y \in {\bf R}^n$, and
\begin{equation}
	T(r \, x) = r \, T(x)
\end{equation}
for all $r \in {\bf R}$ and $x \in {\bf R}^n$.  In particular,
\begin{equation}
	T(0) = 0.
\end{equation}
The identity mapping, which takes $x$ to itself for each $x \in {\bf
R}^n$, is denoted $I$.

	If $T_1$, $T_2$ are linear mappings on ${\bf R}^n$, and $r_1$,
$r_2$ are real numbers, then the linear combination $r_1 \, T_1 + r_2
\, T_2$ is defined as a linear mapping on ${\bf R}^n$, which sends $x
\in {\bf R}^n$ to $r_1 \, T_1(x) + r_2 \, T_2(x)$.  The composition
of $T_1$, $T_2$ is denoted $T_1 \circ T_2$ and defined by
\begin{equation}
	(T_1 \circ T_2)(x) = T_1(T_2(x))
\end{equation}
for $x \in {\bf R}^n$, and is also a linear mapping on ${\bf R}^n$.
Of course if $T$ is a linear mapping on ${\bf R}^n$, then
\begin{equation}
	T \circ I = I \circ T = T.
\end{equation}

	Associated to a linear mapping $T$ on ${\bf R}^n$ we have a $n
\times n$ matrix $(t_{j,k})$ of real numbers, given by
\begin{equation}
	t_{j, k} = \langle T(e_k), e_j \rangle,
\end{equation}
where $e_1, \ldots, e_n$ is the standard basis for ${\bf R}^n$.  This
is equivalent to saying that
\begin{equation}
	T(e_k) = \sum_{j=1}^n t_{j, k} \, e_j.
\end{equation}
Conversely, if $(t_{j, k})$ is an $n \times n$ matrix of real numbers,
then there is a unique linear mapping $T$ on ${\bf R}^n$ associated to
this matrix in this manner.

	The matrix associated to the identity operator $I$ is the
Kronecker delta matrix $(\delta_{j, k})$,
\begin{eqnarray}
	\delta_{j, k} & = & 1 \quad\hbox{when } j = k	\\
		      & = & 0 \quad\hbox{when } j \ne k.
					\nonumber
\end{eqnarray}
If $T_1$, $T_2$ are linear mappings on ${\bf R}^n$ and $(t^1_{j, k})$
and $(t^2_{j, k})$ are the $n \times n$ matrices associated to $T_1$,
$T_2$, respectively, then the matrix associated to $r_1 \, T_1 + r_2
\, T_2$ is given by $(r_1 \, t^1_{j, k} + r_2 \, t^2_{j, k})$.  The $n
\times n$ matrix $(t^3_{j,m})$ associated to the composition $T_3 =
T_1 \circ T_2$ is given by
\begin{equation}
	t^3_{j,m} = \sum_{l=1^n} t^1_{j,l} \, t^2_{l,m},
\end{equation}
which is the usual notion of ``matrix multiplication''.

	If $L$ is any linear subspace of ${\bf R}^n$, then the
\emph{orthogonal complement}\index{orthogonal complement of a linear
subspace of ${\bf R}^n$ in ${\bf R}^n$} of $L$ in ${\bf R}^n$
is denoted $L^\perp$ and defined by
\begin{equation}
	L^\perp = \{w \in {\bf R}^n : w \perp v 
				\hbox{ for all } v \in L\}.
\end{equation}
It is easy to see that $L^\perp$ is also a linear subspace of
${\bf R}^n$.  Furthermore,
\begin{equation}
	L \cap L^\perp = \{0\}.
\end{equation}

	Suppose that $L$ is a linear subspace of ${\bf R}^n$ and that
$v_1, \ldots, v_m$ is an orthonormal basis for $L$.  Define a linear
mapping $P_L$ from ${\bf R}^n$ to itself by
\begin{equation}
	P_L(x) = \sum_{j=1}^m \langle x, v_j \rangle \, v_j.
\end{equation}
This is called the \emph{orthogonal projection of ${\bf R}^n$ onto
$L$}.\index{orthogonal projection of ${\bf R}^n$ onto a linear
subspace}

	One can check that
\begin{equation}
	P_L(x) \in L \quad\hbox{and}\quad x - P_L(x) \in L^\perp
\end{equation}
for each $x \in {\bf R}^n$.  Also, if $y_1$ and $y_2$ are two vectors
in $L$ such that $x - y_1, x - y_2 \in L^\perp$, then $y_1 = y_2$.
In other words, $P_L(x)$ is uniquely determined by the properties
mentioned above, and in particular $P_L(x)$ does not depend on the 
choice of the orthonormal basis $v_1, \ldots, v_m$ for $L$.

	If $w_1, \ldots, w_k$ is an orthonormal collection of vectors
in ${\bf R}^n$ and $u$ is another vector in ${\bf R}^n$ which does
not lie in $\span (w_1, \ldots, w_k)$, then one can use the orthogonal
projection onto $\span (w_1, \ldots, w_k)$ as just described to define
a vector $w_{k+1}$ in ${\bf R}^n$ such that
\begin{equation}
	w_1, \ldots, w_k, w_{k+1}
\end{equation}
is an orthonormal collection of vectors in ${\bf R}^n$ and
\begin{equation}
	\span (w_1, \ldots, w_k, u) 
		= \span (w_1, \ldots, w_k, w_{k+1}).
\end{equation}
Using this one can show that every linear subspace of ${\bf R}^n$ has
an orthonormal basis.  As a result, for each linear subspace of ${\bf
R}^n$ there is an orthogonal projection of ${\bf R}^n$ onto that
subspace.

	If $L$ is a linear subspace of ${\bf R}^n$, then
\begin{equation}
	L \cap L^\perp = \{0\} \quad\hbox{and}\quad 
				L + L^\perp = {\bf R}^n,
\end{equation}
and in particular the sum of the dimensions of $L$ and $L^\perp$
is equal to $n$.  Let us also observe that
\begin{equation}
	(L^\perp)^\perp = L.
\end{equation}
Indeed, $L \subseteq (L^\perp)^\perp$ follows directly from the
definition of the orthogonal complement, and conversely if $x \in
(L^\perp)^\perp$, then $x - P_L(x)$ lies in both $L^\perp$ and
$(L^\perp)^\perp$, and hence is equal to $0$, so that $x = P_L(x) \in
L$.

	By a \emph{bilinear form}\index{bilinear form on ${\bf R}^n$}
we mean a function $B(x, y)$ on ${\bf R}^n \times {\bf R}^n$ with
values in the real numbers which is linear in each of $x$ and $y$,
i.e.,
\begin{equation}
	B(r \, x + r' \, x', y) = r \, B(x, y) + r' \, B(x', y)
\end{equation}
for all $r, r' \in {\bf R}$ and all $x, x', y \in {\bf R}^n$, and
\begin{equation}
	B(x, s \, y + s' \, y') = s \, B(x, y) + s' \, B(x, y')
\end{equation}
for all $s, s' \in {\bf R}$ and all $x, y, y' \in {\bf R}^n$.
If $T$ is a linear mapping on ${\bf R}^n$, then
\begin{equation}
	B(x, y) = \langle T(x), y \rangle
\end{equation}
defines a bilinear form on ${\bf R}^n$.  Conversely, one can check
that every bilinear form on ${\bf R}^n$ arises in this manner.

	If $B(x,y)$ is a bilinear form on ${\bf R}^n$, then we can get
a new bilinear form $B^*(x,y)$ simply by interchanging the variables,
i.e.,
\begin{equation}
	B^*(x, y) = B(y, x).
\end{equation}
If $T$ is a linear transformation on ${\bf R}^n$, then there is a
unique linear transformation $T^*$ on ${\bf R}^n$, called the
\emph{adjoint}\index{adjoint of a linear transformation on ${\bf
R}^n$} of $T$, such that
\begin{equation}
	\langle T^*(x), y \rangle = \langle x, T(y) \rangle.
\end{equation}
If $(t_{j,k})$, $(t^*_{j,k})$ are the matrices associated to $T$,
$T^*$, respectively, then
\begin{equation}
	t^*_{j, k} = t_{k,j}
\end{equation}
for $j, k = 1, \ldots, n$.

	If $r_1$, $r_2$ are real numbers and $T_1$, $T_2$ are linear
mappings on ${\bf R}^n$, then 
\begin{equation}
	(r_1 \, T_2 + r_2 \, T_2)^* = r_1 \, T_1^* + r_2 \, T_2^*.
\end{equation}
Also,
\begin{equation}
	(T_1 \circ T_2)^* = T_2^* \circ T_1^*.
\end{equation}
For any linear transformation $T$ on ${\bf R}^n$,
\begin{equation}
	(T^*)^* = T.
\end{equation}

	A bilinear form $B(x,y)$ on ${\bf R}^n$ is said to be
\emph{symmetric}\index{symmetric bilinear form on ${\bf R}^n$} if $B^*
= B$, a linear mapping $T$ on ${\bf R}^n$ is said to be
\emph{symmetric}\index{symmetric linear mapping on ${\bf R}^n$} if
$T^* = T$, and an $n \times n$ matrix $(t_{j,k})$ is said to be
\emph{symmetric}\index{symmetric matrix} if $t_{j,k} = t_{k,j}$.
These conditions are related in the obvious way, so that a linear
transformation is symmetric if and only if the corresponding bilinear
form is symmetric, and this holds if and only if the associated matrix
is symmetric.  An $n \times n$ matrix $(t_{j,k})$ is said to be
\emph{diagonal}\index{diagonal matrices} if $t_{j,k} = 0$ when $j \ne
k$, which clearly implies that the matrix is symmetric.

	It is easy to see that the adjoint of the identity mapping is
equal to itself.  If $L$ is a linear subspace of ${\bf R}^n$, then the
orthogonal projection $P_L$ of ${\bf R}^n$ onto $L$ is symmetric.  For
if $x, y \in {\bf R}^n$, then
\begin{equation}
	\langle P_L(x), y \rangle = \langle P_L(x), P_L(y) \rangle
			= \langle x, P_L(y) \rangle.
\end{equation}

	Let $T$ be a linear mapping on ${\bf R}^n$.  The
\emph{kernel}\index{kernel of a linear mapping on ${\bf R}^n$} of $T$
is defined to be the set of $x \in {\bf R}^n$ such that $T(x) = 0$,
and the \emph{image}\index{image of a linear mapping on ${\bf R}^n$}
is defined to be the set of $y \in {\bf R}^n$ such that $y = T(x)$ for
some $x \in {\bf R}^n$.  It is easy to see that the kernel and image
of $T$ are linear subspaces of ${\bf R}^n$.

	It is a simple exercise to show that a linear mapping $T$ on
${\bf R}^n$ is one-to-one, which means that $T(x') = T(x)$ implies $x'
= x$ for $x, x' \in {\bf R}^n$, if and only if the kernel of $T$ is
trivial, in the sense that it is equal to $\{0\}$.  One can always
restrict a linear mapping $T$ to the orthogonal complement of its
kernel, and on this subspace of ${\bf R}^n$ $T$ is one-to-one.  Using
this one can check that if $v_1, \ldots, v_m$ is a basis for the
orthogonal complement of the kernel of $T$, then $T(v_1), \ldots,
T(v_m)$ is a basis for the image of $T$.

	As a result, the dimension of the orthogonal complement of
the kernel of $T$ is equal to the dimension of the image of $T$.
Hence the sum of the dimensions of the kernel of $T$ and the image 
of $T$ is equal to $n$.  This implies in turn that $T$ is one-to-one
if and only if $T$ maps ${\bf R}^n$ onto ${\bf R}^n$, which is to
say that the image of $T$ is all of ${\bf R}^n$.

	Let us note as well that if $y \in {\bf R}^n$ lies in the
kernel of the adjoint $T^*$ of a linear mapping $T$ on ${\bf R}^n$,
so that $T^*(y) = 0$, then $y$ lies in the orthogonal complement
of the image of $T$.  Conversely, if $y \in {\bf R}^n$ is in the
orthogonal complement of the image of $T$, so that
\begin{equation}
	\langle y, T(x) \rangle = 0
\end{equation}
for all $x \in {\bf R}^n$, then
\begin{equation}
	\langle T^*(y), x \rangle = 0
\end{equation}
for all $x \in {\bf R}^n$, which implies that $T^*(y) = 0$.
In short, the kernel of $T^*$ is equal to the orthogonal complement
of the image of $T$, and thus the dimensions of the kernels of $T$
and $T^*$ are equal to each other.

	If $T$ is a linear transformation on ${\bf R}^n$ and $L$ is a
linear subspace of ${\bf R}^n$, then $L$ is said to be
\emph{invariant}\index{invariant linear subspace of a linear mapping
on ${\bf R}^n$} under $T$ if $T(L) \subseteq L$, which is to say that
$T(x) \in L$ for all $x \in L$.  One can check that if $L$ is
invariant under $T$, then the orthogonal complement $L^\perp$ of $L$
is invariant under the adjoint $T^*$ of $T$.  As a special case, if
$T$ is a symmetric linear transformation on ${\bf R}^n$ and $L$ is a
linear subspace of ${\bf R}^n$ which is invariant under $L$, then
$L^\perp$ is also invariant under $T$.

	A linear transformation $T$ on ${\bf R}^n$ is said to be
\emph{invertible}\index{invertible linear transformation on ${\bf
R}^n$} if there is another linear transformation $S$ on ${\bf R}^n$
such that
\begin{equation}
	S \circ T = T \circ S = I.
\end{equation}
In this case $S$ is called the \emph{inverse}\index{inverse of a
linear transformation on ${\bf R}^n$} and is denoted $T^{-1}$.  It is
easy to see that the inverse is unique when it exists, and that if
$T_1$, $T_2$ are invertible linear transformations on ${\bf R}^n$,
then $T_1 \circ T_2$ is also invertible, with
\begin{equation}
	(T_1 \circ T_2)^{-1} = T_2^{-1} \circ T_1^{-1}.
\end{equation}

	If $T$ is a linear mapping on ${\bf R}^n$ which is invertible
simply as a mapping, which is to say that there is a mapping $S$ from
${\bf R}^n$ to itself such that $S \circ T = T \circ S = I$, then it
is easy to see that $S$ is linear, so that $T$ is invertible as a
linear mapping.  Of course $T$ is invertible as a mapping on ${\bf
R}^n$ if it is one-to-one and maps ${\bf R}^n$ onto itself.  For
linear mapping these two conditions are equivalent to each other, as
noted before.

	A linear transformation $T$ on ${\bf R}^n$ is said to be an
\emph{orthogonal transformation}\index{orthogonal transformation on
${\bf R}^n$} if
\begin{equation}
	\langle T(x), T(y) \rangle = \langle x, y \rangle
\end{equation}
for all $x, y \in {\bf R}^n$.  Because of the polarization identity,
this is equivalent to saying that
\begin{equation}
	|T(x)| = |x|
\end{equation}
for all $x \in {\bf R}^n$.  One can check that this is also equivalent
to saying that $T$ is invertible, and that
\begin{equation}
	T^{-1} = T^*.
\end{equation}

	If $T$ is a linear mapping on ${\bf R}^n$, then the
\emph{norm}\index{norm of a linear mapping on ${\bf R}^n$} of $T$
is denoted $\|T\|$ and defined by
\begin{equation}
	\|T\| = \sup \{|T(x)| : x \in {\bf R}^n, |x| = 1\}.
\end{equation}
It is easy to check that the nonnegative real numbers $|T(x)|$, $x \in
{\bf R}^n$, $|x| = 1$, are bounded from above, so that the supremum in
this definition makes sense.  To put it another way, the norm $\|T\|$
of $T$ is the nonnegative real number such that
\begin{equation}
	|T(x)| \le \|T\| \, |x|
\end{equation}
for all $x \in {\bf R}^n$, and which is as small as possible.

	Clearly $\|T\| = 0$ if and only if $T = 0$, and the norm of
the identity operator is equal to $1$.  One can check that
\begin{equation}
	\|r \, T\| = |r| \, \|T\|
\end{equation}
for every real number $r$ and every linear operator $T$ on ${\bf
R}^n$, and that
\begin{equation}
	\|T_1 + T_2\| \le \|T_1\| + \|T_2\|
\end{equation}
for all linear mappings $T_1$, $T_2$ on ${\bf R}^n$.  Moreover,
\begin{equation}
	\|T_1 \circ T_2\| \le \|T_1\| \, \|T_2\|
\end{equation}
for all linear mappings $T_1$, $T_2$ on ${\bf R}^n$.

	An equivalent definition of the norm of a linear
transformation $T$ on ${\bf R}^n$ is
\begin{equation}
	\|T\| = \sup \{|\langle T(x), y \rangle| :
				x, y \in {\bf R}^n, |x| = |y| = 1\}.
\end{equation}
Indeed, the right side is less than or equal to $\|T\|$ because
\begin{equation}
	|\langle T(x), y \rangle| \le \|T\|
\end{equation}
when $x, y \in {\bf R}^n$, $|x| = |y| = 1$, by the Cauchy--Schwarz
inequality.  To get the reverse inequality, one can consider the
case where $y$ is a multiple of $T(x)$.

	From this equivalent definition of the norm it follows that
\begin{equation}
	\|T^*\| = \|T\|
\end{equation}
for any linear mapping on ${\bf R}^n$.  Another basic fact is the
``$C^*$ identity''
\begin{equation}
	\|T^* \circ T\| = \|T\|^2.
\end{equation}
This can be checked using the simple formula
\begin{equation}
	\langle (T^* \circ T)(x), y \rangle
		= \langle T(x), T(y) \rangle.
\end{equation}

	If $L$ is a linear subspace of ${\bf R}^n$ which contains a
nonzero element, and $P_L$ is the orthogonal projection of ${\bf R}^n$ 
onto $L$, then
\begin{equation}
	\|P_L\| = 1.
\end{equation}
Indeed, if $x$, $y$ are vectors in ${\bf R}^n$ such that $x \perp y$,
then
\begin{equation}
	|x + y|^2 = |x|^2 + |y|^2,
\end{equation}
and thus
\begin{equation}
	|P_L(u)|^2 \le |P_L(u)|^2 + |u - P_L(u)|^2 = |u|^2
\end{equation}
for all $u \in {\bf R}^n$.  This implies that $\|P_L\| \le 1$, and the
reverse inequalities hold because $P_L(u) = u$ when $u \in L$.

	If $T$ is an orthogonal linear transformation on ${\bf R}^n$,
then of course
\begin{equation}
	\|T\| = 1.
\end{equation}
More precisely, an invertible linear mapping $T$ on ${\bf R}^n$
is an orthogonal transformation if and only if
\begin{equation}
	\|T\| = \|T^{-1}\| = 1.
\end{equation}
This is easy to check from the definitions.

	A symmetric linear transformation $T$ on ${\bf R}^n$ is said
to be \emph{nonnegative}\index{nonnegative symmetric linear
transformation on ${\bf R}^n$} if
\begin{equation}
	\langle T(x), x \rangle \ge 0
\end{equation}
for all $x \in {\bf R}^n$.  More generally, if $A$, $B$ are symmetric
linear transformations on ${\bf R}^n$, then we say that $A$ is less
than or equal to $B$, written
\begin{equation}
	A \le B,
\end{equation}
if $B - A$ is nonnegative.  For any symmetric linear transformation
$T$ on ${\bf R}^n$, we have that
\begin{equation}
	-\|T\| \, I \le T \le \|T\| \, I.
\end{equation}

	If $T$ is a nonnegative symmetric linear operator on ${\bf
R}^n$, then
\begin{equation}
	|\langle T(x), y \rangle|
		\le \langle T(x), x \rangle^{1/2} 
			\, \langle T(y), y \rangle^{1/2}.
\end{equation}
This can be checked in the same manner as for the Cauchy--Schwarz
inequality.  In particular, $T(x) = 0$ when $\langle T(x), x \rangle =
0$.

	If $T$ is a symmetric linear transformation on ${\bf R}^n$,
then
\begin{equation}
	\langle T(x), x \rangle > 0
\end{equation}
for all $x \in {\bf R}^n$ with $x \ne 0$ if and only if $T$ is
invertible.  Indeed, if this condition holds, then the kernel of $T$
contains only the zero vector, and $T$ is invertible.  Conversely, if
$T$ is invertible, then $T(x) \ne 0$ for all $x \in {\bf R}^n$, and
hence the inner product of $T(x)$ and $x$ is nonzero by the result
mentioned in the previous paragraph.

\subsection{Sequences and series}
\label{subsection on sequences and series}

	Fix a positive integer $n$, and let $\{x_j\}_{j=1}^\infty$ be
a sequence of points in ${\bf R}^n$.  This sequence is said to
\emph{converge}\index{convergence of a sequence of points in ${\bf
R}^n$} to a point $x \in {\bf R}^n$ if for every $\epsilon > 0$ there
is an integer $N = N(\epsilon)$ so that
\begin{equation}
	|x_j - x| < \epsilon \quad\hbox{for all } j \ge N.
\end{equation}
In this case $x$ is said to be the \emph{limit}\index{limit of a
sequence of points in ${\bf R}^n$} of the sequence
$\{x_j\}_{j=1}^\infty$, which one can easily show is unique, and one
writes
\begin{equation}
	\lim_{j \to \infty} x_j = x.
\end{equation}

	Sometimes it will be convenient to consider a sequence
$\{x_j\}_{j=a}^\infty$ with a different starting point $a \in {\bf
Z}$, and the same basic notions and results apply just as well.  The
same definition also applies to sequences of complex numbers, using
the standard identification of ${\bf C}$ with ${\bf R}^2$.  In this
regard let us note that a sequence $\{x_j\}_j$ in ${\bf R}^n$
converges to a point $x \in {\bf R}^n$ if and only if for each integer
$l$ such that $1 \le l \le n$ the sequence of $l$th components of the
$x_j$'s converges to the $l$th component of $x$ as a sequence of real
numbers, and in particular a sequence $\{z_j\}_j$ of complex numbers
converges to a complex number $z$ if and only if the real and
imaginary parts of the $z_j$'s converge to the real and imaginary
parts of $z$, respectively.

	Suppose that $\{x_j\}_{j=1}^\infty$, $\{y_j\}_{j=1}^\infty$
are sequences in ${\bf R}^n$ which converge to $x, y \in {\bf R}^n$,
respectively.  A standard result which is easy to verify states that
$\{x_j + y_j\}_{j=1}^\infty$ converges to $x + y$.  In particular this
applies to sequences of real and complex numbers.

	Now suppose that $\{z_j\}_{j=1}^\infty$,
$\{w_j\}_{j=1}^\infty$ are sequences of complex numbers which converge
to the complex numbers $z$, $w$, respectively.  One can then show that
the sequence of products $\{z_j \, w_j\}_{j=1}^\infty$ converges to
the product $z \, w$ of the limits of the individual sequences.  This
includes products of sequences of real numbers as a special case, and
for a pair of sequences $\{x_j\}_{j=1}^\infty$, $\{y_j\}_{j=1}^\infty$
in ${\bf R}^n$ which converge to $x, y \in {\bf R}^n$, respectively,
one has that the sequence of inner products $\{\langle x_j, y_j
\rangle\}_{j=1}^\infty$ converges to the inner product $\langle x, y
\rangle$ of the limits.

	One way to deal with sequences of products is to consider
separately the cases where one of the sequences is constant, and where
one of the sequences converges to $0$.  If one of the sequences tends
to $0$, then it is sufficient for the other sequence to be bounded in
order for the sequence of products to also tend to $0$.  Let us note
as well that if $\{z_j\}_{j=1}^\infty$ is a sequence of nonzero
complex numbers which converges to a nonzero complex number $z$, then
the sequence of reciprocals $\{1/z_j\}_{j=1}^\infty$ converges to the
reciprocal of the limit, $1/z$.

	Sequences of linear transformations on ${\bf R}^n$ can be
converted into sequences of $n \times n$ matrices, which can then be
reinterpreted as sequences in ${\bf R}^m$ with $m = n^2$, so that
convergence can be defined as before.  If $\{S_j\}_{j=1}^\infty$,
$\{T_j\}_{j=1}^\infty$ are sequences of linear transformations on
${\bf R}^n$ which converge to the linear transformations $S$, $T$,
respectively, then $\{S_j +T_j\}_{j=1}^\infty$ and $\{S_j \circ
T_j\}_{j=1}^\infty$ converge to $S + T$ and $S \circ T$, respectively.
If $\{T_j\}_{j=1}^\infty$ is a sequence of invertible linear
transformations on ${\bf R}^n$ which converges to the invertible
linear transformation $T$, then $\{T_j^{-1}\}_{j=1}^\infty$ converges
to $T^{-1}$, as one can show using ``Cramer's rule'' for expressing
the inverse of a matrix in terms of determinants.

	Suppose that $\{x_j\}_{j=1}^\infty$ is a sequence of real
numbers.  We say that this sequence is \emph{monotone
increasing}\index{monotone increasing sequence of real numbers}
if $x_{j+1} \ge x_j$ for all $j$.  Similarly, $\{x_j\}_{j=1}^\infty$
is said to be \emph{monotone decreasing}\index{monotone decreasing
sequence of real numbers} if $x_{j+1} \le x_j$ for all $j$.

	Let $\{x_j\}_{j=1}^\infty$ be a monotone increasing sequence
of real numbers which is bounded from above, which is to say that
there is a real number $b$ such that $x_j \le b$ for all $j$.  In this
case $\{x_j\}_{j=1}^\infty$ converges to a real number $x$, with
\begin{equation}
	x = \sup \{x_j : j \in {\bf Z}_+\}.
\end{equation}
In the same way, if $\{y_j\}_{j=1}^\infty$ is a monotone decreasing
sequence of real numbers which is bounded from below, then
$\{y_j\}_{j=1}^\infty$ converges to $y = \inf_{j \ge 1} y_j$.

	A sequence $\{x_j\}_{j=1}^\infty$ in ${\bf R}^n$ is said to be
a \emph{Cauchy sequence}\index{Cauchy sequence in ${\bf R}^n$} if
for every $\epsilon > 0$ there is a positive integer $N$ such that
\begin{equation}
	|x_j - x_k| < \epsilon \quad\hbox{for all } j, k \ge N.
\end{equation}
It is easy to see that a convergent sequence is a Cauchy sequence, and
in fact the converse holds too, which is to say that every Cauchy
sequence converges.  To show this it is enough to consider a Cauchy
sequence $\{x_j\}_{j=1}^\infty$ of real numbers, and the limit of this
sequence is the same as the limit of the monotone increasing sequence
$\{x'_j\}_{j=1}^\infty$ defined by $x'_j = \inf_{i \ge j} x_i$.

	Now let $\sum_{j=0}^\infty a_j$ be a series of complex
numbers.  Associated to this series is the sequence of partial sums
\begin{equation}
	\sum_{j=0}^n a_j,
\end{equation}
where $n$ runs through the nonnegative integers.  If the sequence of
partial sums converges to a complex number $A$, then we say that the
series $\sum_{j=0}^\infty a_j$ converges to $A$.

	Suppose that $\sum_{j=0}^\infty a_j$, $\sum_{j=0}^\infty b_j$
are two series of complex numbers which converge to the complex
numbers $A$, $B$, respectively.  If $\alpha$, $\beta$ are complex
numbers, then the series
\begin{equation}
	\sum_{j=0}^\infty \alpha \, a_j + \beta \, b_j
\end{equation}
converges to $\alpha \, A + \beta \, B$.  This follows from the
corresponding statement for sequences, since the partial sums of the
new series are equal to the same kind of linear combinations of the
partial sums of the original series.

	The \emph{Cauchy criterion}\index{Cauchy criterion for
convergence of a series of complex numbers} for convergence of a
series states that a series $\sum_{j=0}^\infty a_j$ of complex numbers
converges if and only if for every $\epsilon > 0$ there is a positive
integer $N$ such that
\begin{equation}
	\biggl| \sum_{j=m}^n a_j \biggr| < \epsilon
\end{equation}
for all $m$, $n$ such that $n \ge m \ge N$.  Indeed, this condition is
equivalent to saying that the sequence of partial sums of the series
$\sum_{j=0}^\infty a_j$ is a Cauchy sequence.  A consequence of the
Cauchy criterion is the \emph{comparison test}\index{comparison test
for convergence of a series of complex numbers} for convergence
of a series, which says that if $\sum_{j=0}^\infty a_j$ is a series
of complex numbers and $\sum_{j = 0}^\infty b_j$ is a series of
nonnegative real numbers which converges such that $|a_j| \le b_j$
for all $j$, then $\sum_{j=0}^\infty a_j$ converges.

	Suppose that $\sum_{j=0}^\infty b_j$ is a series whose terms
are nonnegative real numbers.  Then the series converges if and only
if the sequence of partial sums is bounded from above.  This is
because the sequence of partial sums is a monotone increasing sequence
of real numbers.

	The results mentioned in the preceeding two paragraphs are
quite simple, and also quite useful.  In short, the question of
convergence of a series can often be reduced to finding an upper bound
for some collection of nonnegative real numbers.  This is often an
easier task.

	Let us note a basic necessary condition for a series to
converge.  Namely, if $\sum_{j=0}^\infty a_j$ is a series of complex
numbers which converges, then $\{a_j\}_{j=0}^\infty$ converges to $0$.
This is easy to derive from the Cauchy criterion for convergence.

	A series $\sum_{j=0}^\infty a_n$ \emph{converges
absolutely}\index{absolute convergence of a series of complex numbers}
if $\sum_{j=0}^\infty |a_n|$ converges, which implies convergence of
the original series in the ordinary sense.  Absolute convergence is a
more stable kind of convergence.  For example, if a series converges
absolutely, then one can rearrange the terms and still get a series
which converges, and with the same sum.

	A famous test for convergence of series is called the
\emph{Cauchy Condensation Test},\index{Cauchy Condensation Test for
convergence of series} and it says the following.  If
$\sum_{j=1}^\infty b_j$ is a series of nonnegative real numbers such
that the sequence of $b_j$'s is monotone decreasing, then
$\sum_{j=1}^\infty b_j$ converges if and only if the ``condensed''
series $\sum_{k=0}^\infty 2^k b_{2^k}$ converges.  To prove this, one
shows that the partial sums of each can be bounded in terms of the
other.

	Before describing some examples, let us review some auxiliary
facts.  First, the \emph{binomial theorem}\index{binomial theorem}
states that for each positive integer $n$,
\begin{equation}
	(x + y)^n = \sum_{j=0}^n {n \choose j} \, x^j \, y^{n-j},
\end{equation}
where
\begin{equation}
	{n \choose j} = \frac{n!}{j! \, (n-j)!},
\end{equation}
$m! = 1 \cdot 2 \cdots m$, $0! = 1$.  This can be proved by induction,
with the $n = 1$ case being trivial, and let us note also that the
binomial coefficients ${n \choose j}$ are always integers.

	Now let $a$ be a positive real number, and observe that
\begin{equation}
	(1 + a)^n \ge 1 + n \, a
\end{equation}
for all positive integers $n$, by the binomial theorem.  Thus
\begin{equation}
	\lim_{n \to \infty} (1 + a)^{-n} = 0,
\end{equation}
and hence
\begin{equation}
	\lim_{n \to \infty} z^n = 0
\end{equation}
when $z$ is a complex number with $|z| < 1$.  More generally, if $z$
is a complex number with $|z| < 1$ and $k$ is a positive integer, then
\begin{equation}
	\lim_{n \to \infty} n^k \, z^n = 0.
\end{equation}

	If $z$ is a complex number with $|z| < 1$, then the
\emph{geometric series}\index{geometric series} $\sum_{n=0}^\infty
z^n$ converges, and
\begin{equation}
	\sum_{j=0}^\infty z^j = \frac{1}{1 - z}.
\end{equation}
Indeed,
\begin{equation}
	(1 - z) \sum_{j=0}^n = 1 - z^{n+1}
\end{equation}
for each nonnegative integer $n$, so that
\begin{equation}
	\lim_{n \to \infty} \sum_{j=0}^n z^j
		= \lim_{n \to \infty} \frac{1 - z^{n+1}}{1 - z}
			= \frac{1}{1 - z}.
\end{equation}
Of course the series does not converge when $|z| \ge 1$, since $|z|^n
\ge 1$ then for all $n$, and the terms of the series do not tend to
$0$.

	Now suppose that $p$ is a positive real number, and consider
the series
\begin{equation}
	\sum_{n=1}^\infty \frac{1}{n^p}.
\end{equation}
One can show that this series converges if and only if $p > 1$, using
the Cauchy Condensation Test.  More precisely, this reduces the
question to one for geometric series.

	However, the series
\begin{equation}
	\sum_{n=1}^\infty \frac{(-1)^n}{n^p}
\end{equation}
converges for all positive real numbers.  This can be viewed as a
special case of Leibniz' \emph{alternating series
test},\index{alternating series test} which states that a series of
the form
\begin{equation}
	\sum_{n=1}^\infty (-1)^n \, b_n
\end{equation}
converges if $\{b_n\}_{n=1}^\infty$ is a monotone decreasing sequence
of real numbers which converges to $0$.  One can verify this test
using the Cauchy criterion, by estimating the absolute value of sums
of the form $\sum_{n=p}^q (-1)^n \, b_n$ in terms of the absolute
value of the first term $b_p$.

	More generally, a series of the form
\begin{equation}
	\sum_{n=1}^\infty a_n \, b_n
\end{equation}
converges if the sequence of partial sums $A_p = \sum_{n=1}^p a_n$ of
the $a_n$'s is bounded and if $\{b_n\}_{n=1}^\infty$ is a monotone
decreasing sequence of real numbers which converges to $0$.  Indeed,
if we set $A_0 = 0$ for convenience, then for each positive integer
$p$,
\begin{eqnarray}
	\sum_{n=1}^p a_n \, b_n & = & 
			\sum_{n=1}^p (A_n - A_{n-1}) \, b_n	\\
		& = & \sum_{n=1}^p A_n \, b_n 
				- \sum_{n=1}^p A_{n-1} \, b_n
							\nonumber \\
		& = & \sum_{n=1}^p A_n \, b_n
				- \sum_{n=0}^{p-1} A_n \, b_{n+1}
							\nonumber \\
		& = & \sum_{n=1}^{p-1} A_n \, (b_n - b_{n+1})
			+ A_p \, b_p.			\nonumber
\end{eqnarray}
One can check that $\sum_{n=1}^\infty A_n \, (b_n - b_{n+1})$
converges absolutely, since
\begin{equation}
	\sum_{n=1}^\infty (b_n - b_{n+1})
\end{equation}
is a convergent ``telescoping series'' with nonnegative terms, and
hence that $\sum_{n=1}^\infty a_n \, b_n$ converges.

	For example, suppose that $z$ is a complex number such that
$|z| = 1$ and $z \ne 1$, and consider $a_n = z^n$.  In this case we
have that
\begin{equation}
	\sum_{n=1}^p z^n = \frac{z - z^p}{1 - z},
\end{equation}
which is a bounded sequence of complex numbers as $p$ runs through
the positive integers.  Thus we obtain that
\begin{equation}
	\sum_{n=1}^\infty z^n \, b_n
\end{equation}
converges when $z$ is a complex number such that $|z| = 1$ and $z \ne
1$ and $\{b_n\}_{n=1}^\infty$ is a monotone decreasing sequence of
real numbers which converges to $0$.

\subsection{Quarternions and Clifford algebras}
\label{subsection on quarternions and clifford algebras}

	The quarternions are denoted ${\bf H}$, and a quarternion can
be written as $a + b \, i + c \, j + d \, k$, where $a$, $b$, $c$, $d$
are real numbers, and $i$, $j$, $k$ satisfy the conditions
\begin{equation}
	i^2 = j^2 = k^2 = -1
\end{equation}
and
\begin{equation}
	k = i \, j = - j \, i.
\end{equation}
To be more precise, the quarternions are an algebra over the real
numbers, with coordinatewise addition, the real number $1$ being the
multiplicative identity element for the whole algebra, and real
numbers commuting multiplicatively with other quarternions even if
multiplication in ${\bf H}$ is not commutative in general.
Multiplication in ${\bf H}$ does satisfy the usual associative and
distributive rules, and
\begin{equation}
	i \, k = - k \, i = -k, \quad j \, k = k \, j = -i
\end{equation}
in particular.

	Suppose that $x = a + b \, i + c \, j + d \, k$ is a
quarternion.  The \emph{conjugate}\index{conjugate of a quarternion}
of $x$ is denoted $x^*$ and defined by
\begin{equation}
	x^* = a - b \,i - c \, k - d \,k.
\end{equation}
If $x$, $y$ are quarternions, then $(x^*)^* = x$,
\begin{equation}
	(x + y)^* = x^* + y^*,
\end{equation}
and
\begin{equation}
	(x \, y)^* = y^* \, x^*.
\end{equation}

	The \emph{modulus}\index{modulus of a quarternion} of a
quarternion $x = a + b \, i + c \, j + d \, k$ is denoted $|x|$ and
defined by
\begin{equation}
	|x| = (a^2 + b^2 + c^2 + d^2)^{1/2}.
\end{equation}
One can check that
\begin{equation}
	|x|^2 = x \, x^* = x^* \, x,
\end{equation}
and in particular it follows that a nonzero quarternion $x$ is
invertible, with $x^{-1} = |x|^{-2} \, x^*$.  Moreover,
\begin{equation}
	|x \, y| = |x| \, |y|
\end{equation}
for $x, y \in {\bf H}$.

	If $x = a + b \, i + c \, j + d \, k$ is a quarternion, with
$a, b, c, d \in {\bf R}$, as usual, then $a$ is called the \emph{real
part}\index{real part of a quarternion} of $x$.  Of course the real
part of $x$ is equal to the real part of $x^*$.  One can check that if
$x$, $y$ are quarternions, then the real part of $x \, y$ is equal to
the real part of $y \, x$.

	We can identify a quarternion $x = a + b \, i + c \, j + d \,
k$, $a, b, c, d \in {\bf R}$, with the element $(a, b, c, d)$ of ${\bf
R}^4$.  This identification respects the operations of addition and
``scalar'' multiplication by real numbers on ${\bf H}$ and ${\bf
R}^4$, and the modulus of a quarternion corresponds exactly to the
standard Euclidean norm on ${\bf R}^4$.  If $x$, $y$ are quarternions,
then the real part of $x \, y^*$ is the same as the standard inner
product of the cooresponding elements of ${\bf R}^4$.

	Suppose that $\alpha$, $\beta$ are quarternions such that
$|\alpha| = |\beta| = 1$, and consider the mapping
\begin{equation}
	x \mapsto \alpha \, x \, \beta
\end{equation}
from ${\bf H}$ onto itself.  If we identify ${\bf H}$ with ${\bf R}^4$
as in the previous paragraph, then this cooresponds to an orthogonal
linear transformation on ${\bf R}^4$.  This is analogous to the fact
that if $\theta$ is a complex number such that $|\theta| = 1$,
then the mapping 
\begin{equation}
	z \mapsto \theta \, z
\end{equation}
from ${\bf C}$ onto itself corresponds to an orthogonal linear
transformation on ${\bf R}^2$ using the standard identification of
${\bf C}$ with ${\bf R}^2$.

	A quarternion $w$ is said to be
\emph{imaginary}\index{imaginary quarternion} if its real part is
equal to $0$, which is equivalent to
\begin{equation}
	w^* = - w,
\end{equation}
and to
\begin{equation}
	w^2 = - |w|^2.
\end{equation}
If $\alpha$ is a quarternion and $w$ is an imaginary quarternion, then
$\alpha \, w \, \alpha^*$ is also an imaginary quarternion.  If
$\alpha$ is a quarternion with $|\alpha| = 1$, then the mapping
\begin{equation}
	w \mapsto \alpha \, w \, \alpha^* = \alpha \, w \, \alpha^{-1}
\end{equation}
on imaginary quarternions corresponds, under the natural
identification of imaginary quarternions with elements of ${\bf R}^3$,
to an orthogonal linear transformation on ${\bf R}^3$.

	Now let $n$ be a positive integer, and let $\mathcal{C}(n)$
denote the \emph{Clifford algebra}\index{Clifford algebras} with $n$
generators $e_1, \ldots, e_n$.  To be more precise, $\mathcal{C}(n)$
is an algebra over the real numbers which contains a copy of the real
numbers such that the real number $1$ the multiplicative identity
element for the whole Clifford algebra and of course real numbers
commute multiplicatively with all other elements of the Clifford
algebra.  The generators $e_1, \ldots, e_n$ satisfy the relations
\begin{equation}
	e_j^2 = -1
\end{equation}
for $j = 1, \ldots, n$ and
\begin{equation}
	e_j \, e_k = - e_k \, e_j
\end{equation}
when $j \ne k$.

	As an algebra over the real numbers, $\mathcal{C}(n)$ is a
vector space over the real numbers.  As such it has finite dimension.
The basic reason for this is that any product of $e_j$'s can be
reduced to one in which each $e_l$ appears at most once.

	In fact, the dimension of $\mathcal{C}(n)$ as a vector space
over the real numbers is equal to $2^n$.  For let $I$ be a subset of
the set $\{1, \ldots, n\}$, where $I$ consists of
\begin{equation}
	i_1 < i_2 < \cdots < i_l
\end{equation}
for some $l$, $0 \le l \le n$, and let $e_I$ denote the product
\begin{equation}
	e_{i_1} \, e_{i_2} \cdots e_{i_l},
\end{equation}
where this is interpreted as being equal to $1$ when $I = \emptyset$.
The family of $e_I$'s, where $I$ runs through all subsets of $\{1,
\ldots, n\}$, forms a basis for $\mathcal{C}(n)$, which is to say that
every element of $\mathcal{C}(n)$ can be expressed in a unique manner
as a linear combination of the $e_I$'s.

	One can define $\mathcal{C}(n)$ when $n = 0$ to be just the
real numbers themselves.  When $n = 1$, the Clifford algebra
$\mathcal{C}(n)$ is equivalent to the complex numbers.  When $n = 2$,
the Clifford algebra is equivalent to the quarternions.

	Unlike the real numbers, complex numbers, or quarternions, it
is not the case in general that every nonzero element of
$\mathcal{C}(n)$ has a multiplicative inverse.  However, suppose that
$x \in \mathcal{C}(n)$ is of the form
\begin{equation}
	x = x_0 + \sum_{j=1}^n x_j \, e_j,
\end{equation}
where $x_0, x_1, \ldots, x_n \in {\bf R}$.  If $x \ne 0$, then $x$ is
invertible in $\mathcal{C}(n)$, with
\begin{equation}
	x^{-1} = 
  \frac{x_0 - \sum_{j=1}^n x_j \, e_j}{x_0^2 + \sum_{j=1}^n x_j^2}.
\end{equation}

\section{Harmonic and holomorphic functions, 1}
\label{section on harmonic and holomorphic functions, 1}
\setcounter{equation}{0}

	If $x \in {\bf R}^n$ and $r$ is a positive real number, then
we define $B(x, r)$, the open ball with center $x$ and radius $r$
in ${\bf R}^n$\index{open ball $B(x, r)$ in ${\bf R}^n$} by
\begin{equation}
	B(x, r) = \{y \in {\bf R}^n : |y - x| < r\}.
\end{equation}
The closed ball with center $x$ and radius $r$ is denoted
$\overline{B}(x, r)$\index{closed ball $\overline{B}(x, r)$ in ${\bf
R}^n$} and defined by
\begin{equation}
	\overline{B}(x, r)
		= \{y \in {\bf R}^n : |y - x| \le r\}.
\end{equation}
A subset $E$ of ${\bf R}^n$ is said to be \emph{bounded}\index{bounded
subset of ${\bf R}^n$} if it is contained in a ball.

	A subset $U$ of ${\bf R}^n$ is said to be
\emph{open}\index{open subsets of ${\bf R}^n$} if for every $x \in U$
there is a positive real number $r$ such that
\begin{equation}
	B(x, r) \subseteq U.
\end{equation}
The empty set and ${\bf R}^n$ itself are automatically open subsets of
${\bf R}^n$, and one can check that open balls in ${\bf R}^n$ are open
subsets, using the triangle inequality for the Euclidean distance $|v
- w|$, $v, w \in {\bf R}^n$.  The \emph{interior}\index{interior of a
subset of ${\bf R}^n$} of a subset $A$ of ${\bf R}^n$ is denoted
$A^\circ$ and is defined to be the set of $x \in A$ for which there is
a positive real number $r$ such that $B(x, r) \subseteq A$, and one
can check that the interior of $A$ is always an open subset of ${\bf
R}^n$.

	If $E$ is a subset of ${\bf R}^n$ and $p$ is a point in ${\bf
R}^n$, then $p$ is said to be a \emph{limit point}\index{limit point
of a subset of ${\bf R}^n$} of $E$ if for each positive real number
$r > 0$ there is an $x \in E$ such that $x \ne p$ and
\begin{equation}
	|x - p| < r.
\end{equation}
A subset $E$ of ${\bf R}^n$ is said to be \emph{closed}\index{closed
subsets of ${\bf R}^n$} if every point in ${\bf R}^n$ which is a limit
point of $E$ is also an element of $E$.  This is equivalent to saying
that every sequence $\{x_j\}_{j=1}^\infty$ of points in $E$ which
converges to some point $x \in {\bf R}^n$ has its limit $x$ in $E$.

	The empty set and ${\bf R}^n$ itself are automatically closed
subsets of ${\bf R}^n$, and one can check that a closed ball in ${\bf
R}^n$ is a closed subset of ${\bf R}^n$.  If $E$ is a subset of ${\bf
R}^n$, then the \emph{closure}\index{closure of a subset of ${\bf
R}^n$} of $E$ is denoted $\overline{E}$ and is defined to be the union
of $E$ and the set of limit points of $E$, and it is not difficult to
show that the closure of a subset of ${\bf R}^n$ is always a closed
subset of ${\bf R}^n$.  If $A$ is a subset of ${\bf R}^n$, then $A$ is
open if and only if the complement ${\bf R}^n \backslash A$ of $A$ in
${\bf R}^n$, consisting of the points in ${\bf R}^n$ not in $A$, is a
closed subset of ${\bf R}^n$.

\subsection{Some differential operators}
\label{subsection on some differential operators}

	Let $U$ be a nonempty open subset of ${\bf R}^n$, and let
$f(x)$ be a twice continuously-differentiable real or
complex-valued function on $U$.  Put
\begin{equation}
	\Delta f = \sum_{j=1}^n \frac{\partial^2}{\partial x_j^2} f.
\end{equation}
This is called the \emph{Laplacian}\index{Laplace operator} of $f$.

	Now suppose that $n = 2$, and identify ${\bf R}^2$ with ${\bf
C}$.  Writing $z = x + i \, y$, $x, y \in {\bf R}$, for an element of
${\bf C}$, if $f(z)$ is a continuously-differentiable complex-valued
function on a nonempty open subset $U$ of ${\bf C}$, then we put
\begin{equation}
	\frac{\partial}{\partial z} f 
		= \frac{1}{2} \biggl(\frac{\partial}{\partial x} f
			- i \frac{\partial}{\partial y} f \biggr)
\end{equation}
and
\begin{equation}
	\frac{\partial}{\partial \overline{z}} f
		= \frac{1}{2} \biggl(\frac{\partial}{\partial x} f
			+ i \frac{\partial}{\partial y} f \biggr).
\end{equation}
If $f(z)$ is twice continuously-differentiable on $U$, then we have
that
\begin{equation}
   \frac{\partial}{\partial z} \frac{\partial}{\partial \overline{z}} f
= \frac{\partial}{\partial \overline{z}} \frac{\partial}{\partial z} f
= \frac{1}{4} \Delta f.
\end{equation}

	For general $n$ again, assume that $f(x)$ is a
continuously-differentiable function on a nonempty open subset $U$ of
${\bf R}^n$ with values in the Clifford algebra $\mathcal{C}(n)$.
Define the corresponding \emph{left and right Dirac
operators}\index{left and right Dirac operators} acting on $f$ by
\begin{equation}
	\mathcal{D}_L f 
   = \sum_{j=1}^n e_j \, \frac{\partial}{\partial x_j} \, f
\end{equation}
and
\begin{equation}
	\mathcal{D}_R f
   = \sum_{j=1}^n \frac{\partial}{\partial x_j} \, f \, e_j,
\end{equation}
where $e_1, \ldots, e_n$ are the usual multiplicative basis for
$\mathcal{C}(n)$, i.e., we either multiply the first derivatives of
$f$ by the $e_j$'s on the left or on the right, respectively.
If $f$ is twice continuously-differentiable, then we have that
\begin{equation}
	\mathcal{D}_L^2 f = \mathcal{D}_R^2 f = - \Delta f,
\end{equation}
where the Laplacian of $f$, $\Delta f$, is defined for a
Clifford-valued function in the same way as before, which amounts to
applying the Laplacian for real-valued functions to the components of
$f$.

	Let us note that there are some natural variants of this,
where for instance one considers functions on ${\bf R}^{n+1}$ with
values in $\mathcal{C}(n)$, and one uses $1$ as well as $e_1, \ldots,
e_n$ for defining the first-order differential operators.  This type
of set-up is more directly analogous to the one for complex-valued
functions on ${\bf C}$.  For simplicity we shall use the version
described in the previous paragraph.

	If $p$ is a point in ${\bf R}^n$, $r$ is a positive real
number, and $f(x)$ is a continuous function on the closed ball
$\overline{B}(p, r)$ in ${\bf R}^n$, then we write
$\average_{\overline{B}(p, r)} f$ for the average of $f$ over
$\overline{B}(p, r)$, which is the integral of $f$ over
$\overline{B}(p, r)$ divided by the volume of $\overline{B}(p, r)$.
For $j, k = 1, \ldots, n$, we have that
\begin{equation}
	{\average}_{\overline{B}(p, r)} (x_j - p_j)^2
	 = {\average}_{\overline{B}(p, r)} (x_k - p_k)^2,
\end{equation}
and thus
\begin{equation}
	{\average}_{\overline{B}(p, r)} (x_j - p_j)^2
	 = \frac{1}{n} {\average}_{\overline{B}(p, r)} |x - p|^2.
\end{equation}
Also,
\begin{equation}
   {\average}_{\overline{B}(p, r)} |x - p|^2 = \frac{n}{n+2} \, r^2,
\end{equation}
as one can see by reducing to the case where $p = 0$ and using
polar coordinates.

	Let $U$ be a nonempty open subset of ${\bf R}^n$.  Suppose
that $f(x)$ is a twice continuously-differentiable function on $U$,
and that $p$ is an element of $U$.  One can express the Laplacian of
$f$ at $p$ as
\begin{equation}
	\Delta f(p) 
   = \lim_{r \to 0} \frac{n+2}{r^2} 
	\biggl({\average}_{\overline{B}(p, r)} f - f(p) \biggr).
\end{equation}
This follows from the Taylor approximation of $f$ at $p$ of degree
$2$.

	A twice continuously-differentiable function $f$ on a nonempty
open subset $U$ of ${\bf R}^n$ is said to be
\emph{harmonic}\index{harmonic functions on ${\bf R}^n$} if
\begin{equation}
	\Delta f(p) = 0
\end{equation}
for all $p \in U$.  This is equivalent to saying that
\begin{equation}
	\lim_{r \to 0} \frac{1}{r^2} 
	    \biggl({\average}_{\overline{B}(p, r)} f - f(p)\biggr) = 0
\end{equation}
for all $p \in U$.  A stronger ``mean value property'' of harmonic
functions will be described in the next subsection.

	A continuously-differentiable complex-valued function $f(z)$
on an open subset $U \ne \emptyset$ of ${\bf C}$ is said to be
\emph{holomorphic}\index{holomorphic functions on ${\bf C}$} if
\begin{equation}
	\frac{\partial}{\partial \overline{z}} f = 0
\end{equation}
at every point in $U$.  The derivative $f'(z)$ of a holomorphic
function $f(z)$ on $U$ is defined by
\begin{equation}
	f'(z) = \frac{\partial}{\partial z} f(z).
\end{equation}
If we identify ${\bf C}$ with ${\bf R}^2$ in the usual way, then the
differential of a continuously-differentiable complex-valued function
$f(z)$ on $U$ at a point $p$ in $U$ is a linear mapping from ${\bf
R}^2$ to itself, and in fact $(\partial / \partial \overline{z}) f$
vanishes at $p$ if and only if the differential of $f$ at $p$ is given
by multiplication by a complex number on ${\bf C}$, and that complex
number is equal to $(\partial / \partial z) f$ at $p$.

	A continuously-differentiable $\mathcal{C}(n)$-valued function
$f(x)$ on a nonempty open subset $U$ of ${\bf R}^n$ is said to be
\emph{left} or \emph{right Clifford holomorphic},\index{Clifford
holomorphic functions on ${\bf R}^n$} respectively, if
\begin{equation}
	\mathcal{D}_L f(x) = 0
\end{equation}
or
\begin{equation}
	\mathcal{D}_R f(x) = 0
\end{equation}
for all $x \in U$, respectively.  These two conditions are not
equivalent in general, because of noncommutativity, although there are
functions which are both left and right Clifford holomorphic.  These
conditions do not in general imply that the differential of $f$ is
given in terms of Clifford multiplication, which turns out to be a
very restrictive condition.

	If a complex or Clifford holormorphic function is twice
continuously differentiable, then it is also harmonic, in the sense
that its components are harmonic, as one can see from the identities
relating the Laplacian to the various first-order differential
operators mentioned before.  We shall see in the next subsection that the
condition of twice continuous-differentiability actually holds
automatically for complex and Clifford holomorphic functions.  Note
that when $n = 1$, a continuously-differentiable function on an open
interval $(a, b)$ is ``Clifford holomorphic'' if and only if its
derivative is equal to $0$ on the interval, so that the function is
constant, and a twice continuously-differentiable function on $(a, b)$
is harmonic if and only if its second derivative is equal to $0$ on
the interval, which is to say that the function is of the form $\alpha
\, x + \beta$.

	If $h(z)$ is a twice continuously-differentiable function
on a nonempty open subset $U$ of ${\bf C}$ which is harmonic, then
\begin{equation}
	f = \frac{\partial}{\partial z} h
\end{equation}
is a continuously-differentiable function on $U$ which is holomorphic.
If $h(x)$ is a twice continuously-differentiable
$\mathcal{C}(n)$-valued function on a nonempty open subset $U$ of
${\bf R}^n$ which is harmonic, then the functions
\begin{equation}
	f_L = \mathcal{D}_L h \quad\hbox{and}\quad 
		f_R = \mathcal{D}_R h
\end{equation}
are continuously-differentiable $\mathcal{C}(n)$-valued functions on
$U$ which are left and right Clifford holomorphic, respectively.
If $h(x)$ happens to be real-valued, then $f_L = f_R$.

\subsection{Weak derivatives, 1}
\label{subsection on weak derivatives, 1}

	Fix a positive integer $n$, and let $U$ be a nonempty open
subset of ${\bf R}^n$.  Also let $\phi$ be a real, complex, or
$\mathcal{C}(n)$-valued function on $U$.  We say that $\phi$ has
\emph{restricted support}\index{restricted support, for a function on
an open subset of ${\bf R}^n$} in $U$ if the set of $x \in U$ such
that $\phi(x) \ne 0$ is bounded, and if there is a positive real
number $\eta$ such that
\begin{equation}
	|x - y| \ge \eta
\end{equation}
when $x \in U$, $y \in {\bf R}^n$, and $\phi(x) \ne 0$.

	The sum of two functions with restricted support in $U$ has
rerstricted support in $U$, and the product of two functions on $U$,
where at least one of the two functions has restricted support in $U$,
also has restricted support in $U$.  We shall mostly be concerned with
continuous functions on $U$ which have restricted support, and which
then have nice properties in general, such as being bounded, because
of basic results from advanced calculus.  Notice that if $\phi(x)$ is
a continuous function on $U$ with restricted support, and if we extend
$\phi(x)$ to all of ${\bf R}^n$ by setting $\phi(x) = 0$ when $x \in
{\bf R}^n \backslash U$, then this extension is a continuous function
on ${\bf R}^n$ with restricted support in ${\bf R}^n$.

	Suppose that $\phi(x)$ be a continuous function on $U$ with
restricted support.  The usual Riemann integral of $\phi(x)$ on $U$ is
then defined, which we can denote
\begin{equation}
	\int_U \phi(x) \, dx.
\end{equation}
As in the previous paragraph, we can just as well think of $\phi(x)$
as a continuous function on ${\bf R}^n$ with restricted support, and
this integral as an integral over ${\bf R}^n$, even if only a bounded
region in ${\bf R}^n$ is really involved in the integral.

	Let $h(x)$ be a continuous function on $U$.  We say that $h$
is \emph{weakly harmonic}\index{weakly harmonic functions on open
subsets of ${\bf R}^n$} if
\begin{equation}
	\int_U (\Delta \phi(x)) \, h(x) \, dx = 0
\end{equation}
for all twice continuously-differentiable functions $\phi$ on $U$ with
restricted support in $U$.  Notice that if $h$ is twice
continuously-differentiable and harmonic in the usual sense, then $h$
is harmonic in the weak sense, by integration by parts.

	Let us consider for a moment the special case of a function
$\phi(x)$ on ${\bf R}^n$ which is of the form $\rho(|x|^2)$, where
$\rho(t)$ is a twice continuously-differentiable function on the real
line.  By the chain rule,
\begin{equation}
	\Delta \phi(x) 
    = 4 \, |x|^2 \, \rho''(|x|^2) + 2 \, n \, \rho'(|x|^2).
\end{equation}
In other words, $\Delta \phi$ is of the form $\sigma(|x|^2)$, where
\begin{equation}
	\sigma(t) = 4 \, t \, \rho''(t) + 2 \, n \, \rho'(t).
\end{equation}

	Of course we are interested in functions with restricted
support, and so let us assume that there is a positive real number $r$
such that $\rho(t) = 0$ when $t \ge r$.  For simplicity let us assume
also that there is a positive real number $\epsilon$ such that
$\rho(t)$ is constant for $t \le \epsilon$.  In particular,
$\rho'(t)$, $\rho''(t)$, and $\sigma(t)$ are then equal to $0$ for $t
\le \epsilon$.

	By integration by parts, we automatically have
\begin{equation}
	\int_{{\bf R}^n} \Delta \phi(x) \, dx = 0,
\end{equation}
which is equivalent to
\begin{equation}
	\int_0^\infty \sigma(t^2) \, t^{n-1} \, dt = 0,
\end{equation}
and to
\begin{equation}
	\int_0^\infty \sigma(t) \, t^{(n/2) - 1} \, dt = 0.
\end{equation}
Define a function $\theta(u)$ on ${\bf R}$ by
\begin{equation}
	\theta(u) 
   = \frac{1}{4} u^{-n/2} \int_0^u \sigma(t) \, t^{(n/2) - 1} \, dt.
\end{equation}
Thus $\theta(u) = 0$ when $u \le \epsilon$, and also when $u \ge r$.

	By construction,
\begin{equation}
	4 \, u \, \theta'(u) + 2 \, n \, \theta(u) = \sigma(u).
\end{equation}
This is the same equation that $\rho'$ satisfies.  It follows that
$\theta(u) = \rho'(u)$ for all $u \in {\bf R}$.

	Conversely, suppose that $\sigma(t)$ is a continuous function
on the real line which satisifies the conditions mentioned above,
i.e., $\sigma(t) = 0$ when $t \le \epsilon$ and when $t \ge r$, where
$\epsilon$, $r$ are positive real numbers, and the same integral as
before is equal to $0$.  Then we can define $\theta$ again using the
previous formula, so that $\theta$ is a continuous function on ${\bf
R}$ such that $\theta(u) = 0$ when $u \le \epsilon$ and when $t \ge r$
in particular.  We can define $\rho$ on the real line so that $\rho' =
\theta$ and $\rho(u) = 0$ when $u \ge r$, and we also have that
$\rho(u)$ is constant on the set of $u \in {\bf R}$ such that $u \le
\epsilon$.

	In other words, we can reverse the process and start with a
continuous function $\sigma$ on the real line satisfying the
conditions in the previous paragraph, and obtain a twice
continuously-differentiable function $\rho$ on the real line such that
$\rho(u)$ is constant when $u \le \epsilon$ and equal to $0$ when $u
\ge r$.  From the function $\rho$ on the real line we get a function
$\phi(x) = \rho(|x|^2)$ on ${\bf R}^n$.  We want to use these
functions as test functions for the weak harmonicity property.

	Actually, we would like to use translates of such a function
$\phi$.  For each $p \in {\bf R}^n$, put
\begin{equation}
	\phi_p(x) = \phi(x - p) = \rho(|x - p|^2).
\end{equation}
Thus
\begin{equation}
	\Delta \phi_p(x) = \sigma(|x - p|^2).
\end{equation}

	Let us digress a moment with a useful definition.  If $A$ is
a nonempty subset of ${\bf R}^n$ and $x$ is a point in ${\bf R}^n$,
then the \emph{distance from $x$ to $A$}\index{distance from a point
to a nonempty subset of ${\bf R}^n$} is denoted $\dist(x, A)$
and defined by
\begin{equation}
	\dist(x, A) = \inf \{|x - a| : a \in A\}.
\end{equation}
One can check that $\dist(x, A) = 0$ if and only if $x$ lies in the
closure of $A$, and in particular this holds if and only if $x \in A$
when $A$ is a closed subset of ${\bf R}^n$.

	Now suppose that $p$ is an element of $U$, and that $r$ is a
positive real number which satisfies
\begin{equation}
	r^2 < \dist(p, {\bf R}^n \backslash U)
\end{equation}
if $U$ is a proper subset of ${\bf R}^n$, and otherwise is arbitrary
if $U = {\bf R}^n$.  Suppose that $\sigma(t)$ is a continuous function
on the real line such that $\sigma(t) = 0$ when $t \le \epsilon$ for
some $\epsilon > 0$, $\sigma(t) = 0$ when $t \ge r$ for the choice of
$r$ being used now, and $\sigma(t)$ satisfies the integral $0$
condition discussed before.  Thus we get associated functions
$\theta$, $\rho$ on the real line and $\phi_p(x) = \rho(|x - p|^2)$
on ${\bf R}^n$, and $\phi_p$ has restricted support in $U$ because
of the constraint on $r$.

	Let us apply this to our weakly harmonic function $h(x)$ on
$U$.  Namely, we have that
\begin{equation}
	\int_U (\Delta \phi_p(x)) \, h(x) \, dx = 0,
\end{equation}
which is to say that
\begin{equation}
	\int_U \sigma(|x - p|^2) \, h(x) \, dx = 0.
\end{equation}
We want to use this with interesting choices of $\sigma$.

	We can reformulate this condition by saying that if $\psi(x)$
is a continuous radial function on ${\bf R}^n$, so that $\psi(x)$ is
actually a function of $|x|$, $a$ is a positive real number such that
$\psi(x) = 0$ when $|x| \ge a$, $p$ is a point in $U$ such that
\begin{equation}
	\dist(p, {\bf R}^n \backslash U) > a,
\end{equation}
and $\psi$ satisfies the integral condition
\begin{equation}
	\int_{{\bf R}^n} \psi(x) \, dx = 0,
\end{equation}
then
\begin{equation}
	\int_U \psi(x - p) \, h(x) \, dx = 0.
\end{equation}
In other words, $\psi(x)$ corresponds to $\sigma(|x|^2)$.  We have
dropped the condition that $\psi(x)$ be equal to $0$ when $|x|$ is
sufficiently small, because one can reduce to this case through an
approximation argument.

	Let us reformulate this again as follows.  Suppose that
$b_1(x)$, $b_2(x)$ are continuous radial functions on ${\bf R}^n$, $a$
is a positive real number such that $b_i(x) = 0$ when $|x| \ge a$, $i
= 1, 2$, $p$ is an element of ${\bf R}^n$ such that the distance from
$p$ to ${\bf R}^n \backslash U$ is greater than $a$, and
\begin{equation}
	\int_{{\bf R}^n} b_1(x) \, dx 
		= \int_{{\bf R}^n} b_2(x) \, dx = 1.
\end{equation}
Then 
\begin{equation}
	\int_U b_1(x - p) \, h(x) \, dx 
		= \int_U b_2(x - p) \, h(x) \, dx,
\end{equation}
which follows from the statement in the previous paragraph by taking
$\psi(x) = b_1(x) - b_2(x)$.

	This implies that if $b(x)$ is a continuous radial function on
${\bf R}^n$, $a$ is a positive real number such that $b(x) = 0$ when
$|x| \ge a$, $p$ is an element of $U$ whose distance to ${\bf R}^n
\backslash U$ is greater than $a$, and
\begin{equation}
	\int_{{\bf R}^n} b(x) \, dx = 1,
\end{equation}
then
\begin{equation}
	\int_U b(x - p) \, h(x) \, dx = h(p).
\end{equation}
To get this from the previous statement, one can take $b_1(x) = b(x)$,
and choose $b_2(x)$ so that it is concentrated as near to $p$ as one
wants.  Because $h$ is continuous, the integral of $b_2$ times $h$ is
approximately equal to $h(p)$, and in the limit we get the desired
formula.

	In fact, if $p$ is a point in $U$ and $t$ is a positive real
number such that
\begin{equation}
	t < \dist(p, {\bf R}^n \backslash U),
\end{equation}
then
\begin{equation}
	{\average}_{\overline{B}(p, t)} h = h(p).
\end{equation}
This is called the \emph{mean value property}\index{mean value
property of a continuous function on an open subset of ${\bf R}^n$}
for the continuous function $h$ on $U$.  It is easy to derive this
from the previous statement by approximating the average of $h$ on a
closed ball in $U$ by integrals against continuous radial functions,
and conversely one can derive the previous assertion for integrals of
$h$ against continuous radial functions from this one about averages
on closed balls.

	To summarize, a continuous function $h(x)$ on a nonempty
open subset of ${\bf R}^n$ which is weakly harmonic also satisfies
the aforementioned mean value property.

\subsection{Weak derivatives, 2}
\label{subsection on weak derivatives, 2}

	Let $U$ be a nonempty open subset of ${\bf R}^n$, and let $h$
be a continuous function on $U$ which satisfies the mean value
property.  Let $b(x)$ be a continuous radial function on ${\bf R}^n$
and $r$ a positive real number such that $b(x) = 0$ when $|x| \ge r$
and
\begin{equation}
	\int_{{\bf R}^n} b(x) \, dx = 1.
\end{equation}
On the set
\begin{equation}
	U_r = \{p \in U : \dist(p, {\bf R}^n \backslash U) > r\},
\end{equation}
which one can check is an open subset of $U$, we have that
\begin{equation}
	h(p) = \int_U b(x - p) \, h(x) \, dx.
\end{equation}

	We may as well take $b(x)$ to be as smooth as we like here,
and so we choose $b(x)$ so that it is continuously differentiable of
all orders.  By differentiating under the integral sign, it follows
that $h(x)$ is in fact continuously differentiable of all orders
on $U_r$.  This works for all $r > 0$, and therefore $h(x)$ is 
continuously differentiable of all orders on $U$.

	Thus a continuous function which is weakly harmonic is also
continuously differentiable of all orders, and twice
continuously-differentiable in particular.  Such a function is
harmonic in the usual sense.  Similarly, a continuous function which
satisfies the mean value property is harmonic.

	Recall that an open subset $U$ of ${\bf R}^n$ is said to be
\emph{connected}\index{connected open subsets of ${\bf R}^n$} if for
every pair of points $x, y \in U$ there is a continuous path in $U$
that goes from $x$ to $y$; equivalently, $U$ is not connected if it
can be expressed as the union of two disjoint nonempty open subsets of
${\bf R}^n$.  If $h$ is a continuous real-valued function on a
nonempty connected open subset $U$ of ${\bf R}^n$ which satisfies
the mean value property, and if $p$ is a point in $U$ such that
\begin{equation}
	h(x) \le h(p)
\end{equation}
for all $x \in U$, then one can check that
\begin{equation}
	h(x) = h(p)
\end{equation}
for all $x \in U$, because $\{x \in U : h(x) < h(p)\}$ is
automatically an open subset of $U$ since $h$ is continuous, while
$\{x \in U : h(x) = h(p)\}$ is a nonempty open subset of $U$ under the
present conditions.  For the same reasons, if $q$ is a point in $U$
such that
\begin{equation}
	h(x) \ge h(q)
\end{equation}
for all $x \in U$, then 
\begin{equation}
	h(x) = h(q)
\end{equation}
for all $x \in U$.

\subsection{Weak derivatives, 3}
\label{subsection on weak derivatives, 3}

	Let $U$ be a nonempty open subset of ${\bf C}$, and let $f(z)$
be a continuous complex-valued function on $U$.  We say that $f(z)$ is
\emph{weakly holomorphic}\index{weakly holomorphic functions on open
subsets of ${\bf C}$} if
\begin{equation}
   \int_U \biggl(\frac{\partial}{\partial \overline{z}} \psi(z)\biggr)
			\, f(z) \, dz = 0
\end{equation}
for all continuously differentiable complex-valued functions $\psi(z)$
on $U$ with restricted support in $U$.  If $f(z)$ is continuously
differentiable and holomorphic, then $f(z)$ is weakly holomorphic,
by integration by parts.

	Suppose that $f(z)$ is weakly holomorphic on $U$.  We can
choose the test function $\psi(z)$ above to be of the form
\begin{equation}
	\psi(z) = \frac{\partial}{\partial z} \phi(z),
\end{equation}
where $\phi(z)$ is a twice continuously differentiable function
with restricted support in $U$, and it follows that
\begin{equation}
	\int_U (\Delta \phi(z)) \, f(z) \, dz = 0,
\end{equation}
so that $f(z)$ is weakly harmonic.  Thus we obtain that $f(z)$ is
continuously differentiable of all orders on $U$, and that $f(z)$ is
holomorphic in the usual sense.

	Now let $U$ be a nonempty open subset of ${\bf R}^n$, and let
$f(x)$ be a continuous $\mathcal{C}(n)$-valued function on $U$.  We
say that $f(x)$ is \emph{weakly left Clifford
holomorphic}\index{weakly left Clifford holomorphic functions on open
subsets of ${\bf R}^n$} if
\begin{equation}
	\int_U \mathcal{D}_R \psi_1(x) \, f(x) \, dx = 0
\end{equation}
for all continuously-differentiable $\mathcal{C}(n)$-valued functions
$\psi_1$ with restricted support in $U$, and we say that $f(x)$ is
\emph{weakly right Clifford holomorphic}\index{weakly right Clifford
holomorphic functions on open subsets of ${\bf R}^n$} if
\begin{equation}
	\int_U f(x) \, \mathcal{D}_L \psi_2(x) \, dx = 0
\end{equation}
for all continuously-differentiable $\mathcal{C}(n)$-valued functions
$\psi_2$ with restricted support in $U$.  Again one can choose
$\psi_1$, $\psi_2$ to be of the form $\mathcal{D}_R \phi$,
$\mathcal{D}_L \phi$ for a twice continuously-differentiable function
$\phi$ with restricted support in $U$ to obtain that a weakly left or
right Clifford holomorphic function is weakly harmonic, and hence
harmonic, and thus continuously differentiable of all orders and
left or right Clifford holomorphic function in the usual sense.

\section{Harmonic and holomorphic functions, 2}
\label{section on harmonic and holomorphic functions, 2}
\setcounter{equation}{0}

\subsection{Poisson kernels, 1}
\label{subsection on Poisson kernels, 1}

	Fix a positive integer $n$.  By a
\emph{multi-index}\index{multi-index} we mean an $n$-tuple 
\begin{equation}
	\alpha = (\alpha_1, \ldots, \alpha_n)
\end{equation}
where each $\alpha_i$ is a nonnegative integer.  In this case the
\emph{degree}\index{degree of a multi-index} of the multi-index
$\alpha$ is denoted $d(\alpha)$ and defined by
\begin{equation}
	d(\alpha) = \alpha_1 + \cdots + \alpha_n.
\end{equation}

	Let $\alpha = (\alpha_1, \ldots, \alpha_n)$ be a multi-index.
The corresponding \emph{monomial}\index{monomials on ${\bf R}^n$} is
the function of $x = (x_1, \ldots, x_n) \in {\bf R}^n$ defined by
\begin{equation}
	x^\alpha = x_1^{\alpha_1} \cdots x_n^{\alpha_n}.
\end{equation}
The \emph{degree}\index{degree of a monomial on ${\bf R}^n$} of the
monomial $x^\alpha$ is defined to be the degree of the multi-index
$\alpha$.

	By a \emph{polynomial}\index{polynomial on ${\bf R}^n$} on
${\bf R}^n$ we mean a function which is a finite linear combination of
monomials.  A \emph{harmonic polynomial}\index{harmonic polynomial on
${\bf R}^n$} is simply a polynomial which is harmonic.  For a
polynomial $p(x)$ on ${\bf R}^n$, the Laplacian $\Delta p(x)$ can
be expressed algeraically, and thus the property of being harmonic
can also be expressed algebraically.

	A polynomial $p(x)$ on ${\bf R}^n$ is said to be
\emph{homogeneous of degree $d$},\index{homogeneous polynomials on
${\bf R}^n$} where $d$ is a nonnegative integer, if $p(x)$ is a finite
linear combination of monomials $x^\alpha$ with $d(\alpha) = d$.  This
is equivalent to saying that
\begin{equation}
	p(t \, x) = t^d \, p(x)
\end{equation}
for all $t \in {\bf R}$ and all $x \in {\bf R}^n$.  Every polynomial
on ${\bf R}^n$ can be decomposed into homogeneous parts in a simple
way.

	If $p(x)$ is a homogeneous polynomial on ${\bf R}^n$ of degree
$d$, then the Lapacian $\Delta p(x)$ of $p(x)$ is equal to $0$ if $d =
0, 1$ and is a homogeneous polynomial of degree $d - 2$ when $d \ge
2$.  If $h(x)$ is a harmonic polynomial on ${\bf R}^n$, then its
homogeneous components are also harmonic.  If $q(x)$ is a homogeneous
polynomial of degree $\ell$, then $|x|^2 \, q(x)$ is a homogeneous
polynomial of degree $\ell + 2$.

	Suppose that $p(x)$ is a homogeneous polynomial on ${\bf R}^n$
of degree $d$.  Then
\begin{equation}
	\sum_{j=1}^n x_j \, \frac{\partial}{\partial x_j} \, p(x)
		= d \, p(x),
\end{equation}
which is known as ``Euler's identity''.  This follows by computing the
derivative of $p(t \, x)$ in $t$ and setting $t = 1$, where the
derivative is computed using the chain rule, and by using homogeneity
to write $p(t, x)$ as $t^d \, p(x)$ and differentiating $t^d$
directly.

	If $q(x)$ is a polynomial on ${\bf R}^n$, then
\begin{equation}
	\Delta (|x|^2 \, q(x)) 
		= 2 \, n \, q(x) 
   + 4 \sum_{j=1}^n x_j \, \frac{\partial}{\partial x_j} \, q(x)
	+ |x|^2 \, \Delta q(x).
\end{equation}
Assuming also that $q(x)$ is homogeneous of degree $\ell$, it follows
that
\begin{equation}
	\Delta(|x|^2 \, q(x)) = 2 \, (n + 2 \, \ell) \, q(x)
					+ |x|^2 \, \Delta q(x).
\end{equation}
Supposing further that $q(x)$ is harmonic, we obtain that
\begin{equation}
	\Delta(|x|^2 \, q(x)) = 2 \, (n + 2 \, \ell) \, q(x),
\end{equation}
and more generally $\Delta(|x|^{2k} \, q(x))$ is a positive integer
multiple of $|x|^{2k-2} \, q(x)$ when $k$ is a positive integer.

	A key result now is that if $p(x)$ is any polynomial on ${\bf
R}^n$, then $p(x)$ can be written as a finite linear combination of
polynomials of the form $|x|^{2k} \, h(x)$, where $k$ is a nonnegative
integer and $h(x)$ is a harmonic polynomial.  To be more precise, if
$p(x)$ is a homogeneous polynomial of degree $d$, then $p(x)$ can be
expressed as a finite linear combination of polynomials of the form
$|x|^{2k} \, h(x)$, where $k$ is a nonnegative integer and $h(x)$ is a
homogeneous harmonic polynomial of degree $\ell$.  This can be
verified using the computations above.

\subsection{Poisson kernels, 2}
\label{subsection on Poisson kernels, 2}

	Fix a positive integer $n$, and let ${\bf S}^{n-1}$ denote the
unit sphere in ${\bf R}^n$, so that
\begin{equation}
	{\bf S}^{n-1} = \{x \in {\bf R}^n : |x| = 1\}.
\end{equation}
Also let $f(x)$ be a continuous real-valued function on ${\bf
S}^{n-1}$.  We are interested in continuous real-valued functions
$h(x)$ on the closed unit ball $\overline{B}(0,1)$ which agree with
$f(x)$ on the boundary ${\bf S}^{n-1}$ and which are harmonic on the
open unit ball $B(0,1)$.

	Suppose that $h(x)$ is a continuous real-valued function on
$\overline{B}(0,1)$ which is harmonic on $B(0,1)$.  By standard
results from advanced calculus, $h(x)$ attains its maximum and minimum
on $\overline{B}(0,1)$, and we have seen that if the maximum or
minimum are attained in the interior $B(0,1)$, then $h(x)$ is
constant.  As a consequence, the maximum and minimum are always
attained on the boundary ${\bf S}^{n-1}$, and it follows that
\begin{equation}
	\sup \{|h(x)| : x \in \overline{B}(0,1)\}
		= \sup \{|h(x)| : x \in {\bf S}^{n-1} \},
\end{equation}
and in particular that $h(x) = 0$ for all $x \in \overline{B}(0,1)$
when $h(x) = 0$ for all $x \in {\bf S}^{n-1}$.

	Thus a continuous real-valued function $h(x)$ on
$\overline{B}(0,1)$ which is harmonic on $B(0,1)$ is uniquely
determined by its boundary values, i.e., if $h_1(x)$, $h_2(x)$ are two
continuous real-valued functions on $\overline{B}(0,1)$ which are
harmonic on $B(0,1)$ and which are equal on ${\bf S}^{n-1}$, then
$h_1(x) - h_2(x)$ is a continuous real-valued function on
$\overline{B}(0,1)$ which is harmonic on $B(0,1)$ and equal to $0$ on
${\bf S}^{n-1}$, and hence is equal to $0$ on all of
$\overline{B}(0,1)$.  Let us now consider the question of existence of
continuous functions on $\overline{B}(0,1)$ which are harmonic on
$B(0,1)$ and whose values on the boundary ${\bf S}^{n-1}$ are
prescribed in advance.  When $n = 1$ the existence result is
immediate, because ${\bf S}^{n-1}$ consists of the two points $1$,
$-1$, and affine functions are harmonic.

	In general, if $f(x)$ is a continuous real-valued function on
${\bf S}^{n-1}$, then there is a sequence of polynomials
$\{p_j(x)\}_{j=1}^\infty$ on ${\bf R}^n$ such that the restrictions of
the $p_j(x)$'s to ${\bf S}^{n-1}$ converge uniformly to $f(x)$.  This
is a well-known result from advanced calculus.  By the results of the
previous subsection, there are harmonic polynomials $h_j(x)$ on ${\bf
R}^n$ which agree with the $p_j(x)$'s on ${\bf S}^{n-1}$.

	For positive integers $j$, $k$ we have that
\begin{eqnarray}
\lefteqn{\sup \{|h_j(x) - h_k(x)| : x \in \overline{B}(0,1)\}} \\
	& & = \sup \{|h_j(x) - h_k(x)| : x \in {\bf S}^{n-1}\}
						\nonumber \\
	& & = \sup \{|p_j(x) - p_k(x)| : x \in {\bf S}^{n-1}\}.
						\nonumber
\end{eqnarray}
Using this one can check that $\{h_j(x)\}_{j=1}^\infty$ converges
uniformly to a continuous function $h(x)$ on $\overline{B}(0,1)$ which
agrees with $f(x)$ on the boundary ${\bf S}^{n-1}$.  Also, $h(x)$ is
harmonic on $B(0,1)$, because it follows from uniform convergence that
$h(x)$ is weakly harmonic.

	One can also look for an integral formula for a harmonic
extension on the closed unit ball of a continuous function on the unit
sphere.  Namely, if $f(x)$ is a continuous function on ${\bf
S}^{n-1}$, one would like to write
\begin{equation}
	h(x) = \int_{{\bf S}^{n-1}} P(x,y) \, f(y) \, dy
\end{equation}
for a harmonic function $h(x)$ on $B(0,1)$ with boundary values
$f(x)$.  Here $dy$ denotes the usual element of surface integration
in ${\bf R}^n$.

	More precisely, one would like $P(x, y)$ to be a continuous
function defined for $x \in B(0 ,1)$ and $y \in {\bf S}^{n-1}$.  In
order for $h(x)$ to be harmonic, $P(x, y)$ should be harmonic as a
function of $x \in B(0,1)$ for each $y \in {\bf S}^{n-1}$.  In order
for the boundary values of $h$ to be given by $f$ on ${\bf S}^{n-1}$,
for each $z \in {\bf S}^{n-1}$ the function $P(x, y)$, as a function
of $y$ on ${\bf S}^{n-1}$, should converge to the Dirac delta function
$\delta_z(y)$ at $z$ as $x \in B(0,1)$ tends to $z$, in a suitable
sense, which means in particular that $P(x, y)$ should tend to $0$
for $y \ne z$ as $x$ tends to $z$.

	Because $h(x)$ should be equal to $1$ for all $x \in B(0,1)$
when $f(y) = 1$ for all $y \in {\bf S}^{n-1}$, we should have
\begin{equation}
	\int_{{\bf S}^{n-1}} P(x, y) \, dy = 1
\end{equation}
for all $x \in B(0, 1)$.  If $f(y)$ is a nonnegative real-valued
function on ${\bf S}^{n-1}$, then $h(x)$ should be a nonnegative
real-valued function on $B(0,1)$, and this leads to the condition that
$P(x, y)$ should be a nonnegative real number for all $x \in B(0, 1)$
and $y \in {\bf S}^{n-1}$.  Also, $P(0, y)$ should be equal to
$1/\nu_{n-1}$ for all $y \in {\bf S}^{n-1}$, where $\nu_{n-1}$
denotes the area of the unit sphere ${\bf S}^{n-1}$.

	In fact there is a function $P(x, y)$ with these properties,
called the \emph{Poisson kernel for the unit ball}.\index{Poisson
kernel for the unit ball} Specifically, $P(x, y)$ is given by the
formula
\begin{equation}
	P(x, y) = 
	   \frac{1}{\nu_{n-1}} \, \frac{1 - |x|^2}{|x - y|^n}.
\end{equation}
The Poisson kernel is unique, and it provides another approach to the
existence of continuous functions on the closed unit ball which are
harmonic on the open unit ball and have prescribed boundary values.

	Similarly, on any closed ball $\overline{B}(w, r)$ a
continuous function on the boundary sphere can be extended to a unique
continuous function on the closed ball which is harmonic on the
corresponding open ball, and the extension on the open ball can be
expressed by an analogous Poisson integral of the boundary values.
Let us note that the Poisson kernel is real-analytic as a function of
$x$, which is to say that it admits local representations in terms of
power series.  Each harmonic function on a nonempty subset of ${\bf
R}^n$ is locally given as a Poisson integral, and as a result it
follows that each harmonic function is also real analytic.

	Now let us describe an application of the existence of
harmonic functions on balls in ${\bf R}^n$ with prescribed continuous
boundary values.  Let $U$ be a nonempty open subset of ${\bf R}^n$,
and let $h(x)$ be a continuous function on $U$.  We say that $h(x)$
satisfies the \emph{weak mean value property}\index{weak mean value
property, for a continuous function on an open subset of ${\bf R}^n$}
if for each $x \in U$ there is a positive real number $r(x)$ such that
the usual mean value property holds for $h$ for balls centered at $x$
and with radius less than $r(x)$, which is to say that
\begin{equation}
	{\average}_{\overline{B}(x, r)} h = h(x)
\end{equation}
for all positive real numbers $r$ such that $r < r(x)$.

	Suppose that $U$ is a nonempty connected open subset of ${\bf
R}^n$, and that $h(x)$ is a continuous real-valued function on $U$
which satisfies the weak mean value function.  If $p$ is a point in
$U$ at which $h$ attains its maximum, so that
\begin{equation}
	h(x) \le h(p)
\end{equation}
for all $x \in U$, then one can show that $h$ is constant on $U$, with
$h(x) = h(p)$ for all $x \in U$.  Namely, $\{x \in U : h(x) < h(p)\}$
is an open subset of $U$ simply because $h$ is continuous, while $\{x
\in U : h(x) = h(p) \}$ is an open subset of $U$ because of the weak
mean value property and the assumption that $h$ attains its maximum at
$p$, and hence $\{x \in U : h(x) = h(p)\} = U$ by the connectedness of
$U$.

	Another feature of the weak mean value property is that it is
preserved by restricting a function to a smaller open set.  As a
result of this and the maximum principle in the preceding paragraph,
it follows that a function which satisfies the weak mean value
property on a nonempty open subset of ${\bf R}^n$ can be represented
locally as a harmonic function on a ball, namely a harmonic function
whose boundary values are the given function.  In short, a continuous
function on a nonempty open subset of ${\bf R}^n$ which satisfies
the weak mean value property is actually harmonic.

\subsection{Poisson kernels, 3}
\label{subsection on Poisson kernels, 3}

	Fix a positive integer $n$ again, and let us identify ${\bf
R}^{n+1}$ with the Cartesian product ${\bf R}^n \times {\bf R}$,
and write $(x, t)$, $x \in {\bf R}^n$, $t \in {\bf R}$, for an
element of ${\bf R}^n$.  The \emph{upper half-space}\index{upper
half-space in ${\bf R}^{n+1}$} in ${\bf R}^{n+1}$ is denoted
$\mathcal{U}^{n+1}$ and defined by
\begin{equation}
	\mathcal{U}^{n+1} = \{ (x, t) \in {\bf R}^n \times {\bf R} :
				t > 0 \},
\end{equation}
and the corresponding closed half-space is given by
\begin{equation}
	\overline{\mathcal{U}}^{n+1}
		= \{ (x, t) \in {\bf R}^n \times {\bf R} : t \ge 0\}.
\end{equation}
We shall be interested in continuous functions on the closed
half-space which are harmonic on the open half-space.

	Let us begin with a version of the \emph{Schwarz reflection
principle}.\index{Schwarz reflection principle} Namely, if $h(x, t)$
is a continuous function on the closed upper half-space which is
harmonic on the open upper half-space, and if
\begin{equation}
	h(x, 0) = 0
\end{equation}
for all $x \in {\bf R}^n$, then we can extend $h(x, t)$ to all of
${\bf R}^n \times {\bf R}$ by setting
\begin{equation}
	h(x, -t) = - h(x, t)
\end{equation}
when $t > 0$, and this extension is harmonic.  Indeed, it is easy to
see that this extension satisfies the weak mean value property on
${\bf R}^n \times {\bf R}$.

	Next, suppose that $f(x)$ is a harmonic function on ${\bf
R}^m$ which is \emph{bounded},\index{bounded harmonic functions on
${\bf R}^n$} so that there is a positive real number $C$ so that
\begin{equation}
	|f(x)| \le C
\end{equation}
for all $x \in {\bf R}^m$.  In this case one can show that $f(x)$ is a
constant function.  Indeed, for any $x, y \in {\bf R}^m$ and any
positive real number $r$, we have that
\begin{equation}
	f(x) - f(y) = {\average}_{\overline{B}(x, r)} f
			- {\average}_{\overline{B}(x, r)} f,
\end{equation}
and one can show that the right side tends to $0$ as $r \to \infty$
because $f$ is bounded.

	In fact, there is a generalization of this for harmonic
functions $f(x)$ on ${\bf R}^m$ which are of at most polynomial
growth, in the sense that there is a positive real number $C$
and a positive integer $\ell$ such that
\begin{equation}
	|f(x)| \le C \, (1 + |x|)^\ell
\end{equation}
for all $x \in {\bf R}^m$.  In this case one can show that $f(x)$ is a
polynomial of degree less than or equal to $\ell$.  For instance, one
can show that the $\ell$th order derivatives of $f(x)$ are bounded
harmonic functions on ${\bf R}^m$, and hence are constant.

	At any rate, we conclude that if $h_1(x, t)$, $h_2(x, t)$ are
bounded continuous functions on the closed upper half-space in ${\bf
R}^m \times {\bf R}$ which are harmonic on the open half-space and
which are equal on the boundary, then they are equal on the whole
closed half-space.  Indeed, under these conditions $h_1(x, t) - h_2(x,
t)$ is a bounded continuous function on the closed half-space which is
harmonic on the open half-space and equal to $0$ when $t = 0$, and one
can use reflection to get a bounded harmonic function on all of ${\bf
R}^n \times {\bf R}$, which is therefore constant, and hence equal to
$0$.  Notice that $h(x, t) = t$ is a harmonic function on ${\bf R}^n
\times {\bf R}$ which is equal to $0$ when $t = 0$ and which is not
equal to $0$ when $t \ne 0$.

	Now let us consider the corresponding question of existence of
bounded continuous functions on the closed half-space which are
harmonic on the open half-space and which have prescribed boundary
values.

	We can begin with the point of view of Fourier analysis.
Namely, suppose that $\xi$ is an element of ${\bf R}^n$, and consider
the bounded continuous function
\begin{equation}
	\exp (i \xi \cdot x)
\end{equation}
on ${\bf R}^n$, as a function of $x$.  One can check that
\begin{equation}
	\exp (- |\xi| t) \exp(i \xi \cdot x)
\end{equation}
is a harmonic function of $(x, t)$ on ${\bf R}^n \times {\bf R}$ which
has modulus less than or equal to $1$ when $t \ge 0$ and which agrees
with the specified function of $x$ when $t = 0$.

	As in the case of the unit ball, there is a nice integral
formula in this case.  Namely, if $f(x)$ is a bounded continuous
function on ${\bf R}^n$, then
\begin{equation}
	h(x, t) = \int_{{\bf R}^n} P(x - y, t) \, f(y) \, dy
\end{equation}
defines a bounded harmonic function on the upper-half space whose
boundary values are equal to $f(x)$, where $P(x, t)$ is the
\emph{Poisson kernel}\index{Poisson kernel for the upper half-space}
for the upper half-space, which is the function on the upper
half-space given by
\begin{equation}
	a_n \, \frac{t}{(|x|^2 + t^2)^{(n+1)/2}},
\end{equation}
with $a_n$ a positive real number.

	To be more precise, if $f(x)$ is a bounded continuous function
on ${\bf R}^n$, and we define $h(x, t)$ on the upper half-space using
the formula above, and set $h(x, 0) = f(x)$, then $h(x, t)$ is a
bounded continuous function on the closed half-space which is harmonic
on the open half-space.  Harmonicity of $h(x, t)$ is a consequence of
harmonicity of $P(x, t)$, as a function on the upper half-space.  As a
function of $x$, $P(x, t)$ tends to the Direc delta function at $0$ as
$t \to 0$ in a suitable sense, and this corresponds to the fact that
$h(w, t)$ tends to $f(w)$ as $t \to 0$.

	As in the case of the unit ball, $P(x, t)$ is positive
everywhere, and this corresponds to the fact that $h(x, t)$ is
positive everywhere when $f(x)$ is nonnegative and not identically
equal to $0$.  Also,
\begin{equation}
	\int_{{\bf R}^n} P(x, t) \, dx = 1
\end{equation}
for all $t > 0$, which corresponds to the fact that $h(x, t) = 1$
for all $(x, t)$ in the upper half-space when $f(x) = 1$ for all
$x \in {\bf R}^n$.  The Poisson kernel $P(x, t)$ is the unique
function on the upper half-space which represents bounded harmonic
functions on the upper half-space that are continuous up to the
boundary in this manner.

\subsection{Fundamental solutions}
\label{subsection on fundamental solutions}

	Fix a positive integer $n$.  Let us first consider the
question of finding a \emph{Newtonian kernel}\index{Newtonian kernel}
$N_n(x)$ in dimension $n$, which should be a function on ${\bf R}^n$
whose Laplacian is equal to the Direc delta function at $0$.  This is
to be interpreted in a weak sense, which is to say that
\begin{equation}
	\int_{{\bf R}^n} N_n(x) \, \Delta \phi(x) \, dx 
			= \phi(0)
\end{equation}
when $\phi(x)$ is a twice continuously-differentiable function
on ${\bf R}^n$ with restricted support.

	On ${\bf R}^n \backslash \{0\}$, $N_n(x)$ should be a harmonic
function.  This does not mean that $N_n(x)$ is necessarily continuous
at $0$, however.  It turns out that the singularity of $N_n(x)$ at the
origin is mild enough so that the above integral makes sense as an
principal value integral, and behaves reasonably well, with absolute
convergence in particular.

	The Newtonian kernel is not uniquely determined by the
condition that its Laplacian be equal to the Dirac delta function at
$0$ in the weak sense, because one could add to $N_n(x)$ any harmonic
function on ${\bf R}^n$ and get a function with the same property.
However, if one imposes the condition that $N_n(x)$ be a radial
function, then it is unique except for adding a constant.  This is
because radial harmonic functions on ${\bf R}^n \backslash \{0\}$
are characterized by a homogeneous second-order ordinary differential
equation, with a two-dimensional space of solutions, and constant
functions account for a one-dimensional space of solutions.

	When $n = 1$, one can take
\begin{equation}
	N_1(x) = \frac{1}{2} \, |x|.
\end{equation}
When $n = 2$ one can take
\begin{equation}
	N_2(x) = b_2 \log|x|,
\end{equation}
and when $n \ge 3$ one can take
\begin{equation}
	N_n(x) = - b_n |x|^{2 - n},
\end{equation}
where $b_n$, $n \ge 2$, are positive real numbers.  This is not too
difficult to check.

	If $y$ is a point in ${\bf R}^n$, then $N_n(x - y)$ has
Laplacian, as a function of $x$, equal to the Dirac delta function at
$y$ in the weak sense, so that
\begin{equation}
	\int_{{\bf R}^n} N_n(x - y) \, \Delta \phi(x) \, dx
					= \phi(y)
\end{equation}
for all twice continuously-differentiable functions $\phi$ on ${\bf
R}^n$ with restricted support.  Another way to look at this is that
if $\psi$ is a continuous function on ${\bf R}^n$ with restricted
support, then
\begin{equation}
	\Psi(x) = \int_{{\bf R}^n} N_n(x - y) \, \psi(y) \, dy
\end{equation}
is a continuous function on ${\bf R}^n$ such that
\begin{equation}
	\Delta \Psi = \psi
\end{equation}
in the weak sense, which is to say that
\begin{equation}
	\int_{{\bf R}^n} \Psi(x) \, \Delta \phi(x) \, dx
	   = \int_{{\bf R}^n} \psi(x) \, \phi(x) \, dx
\end{equation}
for all twice continuously-differentiable functions $\phi$ on ${\bf
R}^n$ with restriced support.  In general, $\Psi$ is not quite twice
continuously-differentiable when $\psi$ is continuous, but this is
almost true, and $\Psi$ is twice continuously-differentiable under
additional mild conditions on $\psi$.

	Of course the idea of a fundamental solution like this
makes sense for differential operators in general.

	For instance, when $n = 1$, one can consider the derivative of
$N_1(x)$, which is equal to $-1/2$ when $x < 0$ and to $1/2$ when $x >
0$.  One can make the convention of setting this function to be $0$ at
$x = 0$.  At any rate, $N_1'(x)$ is a fundamental solution for
$d/dx$ in the sense that its derivative is equal to the Dirac delta
function at $0$ in the weak sense, so that
\begin{equation}
	- \int_{\bf R} N_1'(x) \, \phi'(x) \, dx
			= \phi(0)
\end{equation}
for all continuously-differentiable functions on the real line with
restricted support.

	Similarly, a nonzero multiple of $1/z$ defines a fundamental
solution for the differential operator $\partial / \partial
\overline{z}$.  In other words, the $\partial / \partial \overline{z}$
derivative of $1 / z$ is equal to a nonzero constant times the Dirac
delta function at $0$ on ${\bf C}$, in the sense that
\begin{equation}
	\int_{\bf C} \frac{1}{z} \, 
		\frac{\partial}{\partial \overline{z}} \phi(z) \, dz
\end{equation}
is equal to a nonzero constant times $\phi(0)$ for all
continuously-differentiable functions $\phi$ on ${\bf C}$ with
restricted support.  Note that $1/z$ is indeed holomorphic on ${\bf C}
\backslash \{0\}$, as it should be to have $\partial / \partial
\overline{z}$ derivative equal to $0$ there, and that $1/z$ is a
constant multiple of $\partial / \partial z N_2(z)$, identifying
${\bf R}^2$ with ${\bf C}$ in the usual manner.

	For general $n$, the Clifford-valued function
\begin{equation}
	E_n(x) = \frac{\sum_{j=1}^n x_i \, e_i}{|x|^n}
\end{equation}
is a nonzero constant multiple of a fundamental solution for both
$\mathcal{D}_L$ and $\mathcal{D}_R$.  This means that
\begin{equation}
	\int_{{\bf R}^n} \mathcal{D}_R \phi(x) \, E_n(x) \, dx
\end{equation}
and
\begin{equation}
	\int_{{\bf R}^n} E_n(x) \, \mathcal{D}_L \phi(x) \, dx
\end{equation}
are equal to a constant multiple of $\phi(0)$ whenever $\phi(x)$ is a
continuously-differentiable Clifford-valued function on ${\bf R}^n$
with restricted support.  This function $E_n(x)$ is a constant
multiple of $\mathcal{D}_L N_n(x) = \mathcal{D}_R N_n(x)$, and is both
left and right Clifford holomorphic on ${\bf R}^n \backslash \{0\}$,
as it should be.

\section{Harmonic and holomorphic functions, 3}
\label{section on harmonic and holomorphic functions, 3}
\setcounter{equation}{0}

\subsection{Subharmonic functions, 1}
\label{subsection on subharmonic functions, 1}

	In dimension $1$, the notion of subharmonicity reduces to
that of convexity.  Let $(a, b)$ be an open interval in the real
line, which is to say the set of real numbers $x$ such that
\begin{equation}
	a < x < b,
\end{equation}
where of course we assume that $a < b$, and also we allow $a$ to
be $-\infty$ or a finite real number, and $b$ to be a finite real
number or $+\infty$, so that our interval may be unbounded, like the
set of positive real numbers, or the whole real line itself.
A real-valued function $f(x)$ is said to be convex if
\begin{equation}
	f(\lambda \, x + (1 - \lambda) \, y)
		\le \lambda \, f(x) + (1 - \lambda) \, f(y)
\end{equation}
for all $x, y \in (a, b)$ and all $\lambda \in (0,1)$.

	It is a well-known exercise that convex functions are
continuous.  In the other direction, if one assumes that $f(x)$ is a
continuous real-valued function on $(a, b)$, then more restrictive
convexity conditions imply the general one, e.g., it is enough to take
$\lambda = 1/2$.  If $f(x)$ is a twice differentiable function on $(a,
b)$, then $f(x)$ is convex if and only if $f''(x) \ge 0$ for all $x
\in (a, b)$.

	For any continuous real-valued function $f(x)$ on $(a, b)$,
we can say that $f'' \ge 0$ on $(a, b)$ in the weak sense if
\begin{equation}
	\int_a^b f(x) \, \phi''(x) \, dx \ge 0
\end{equation}
for all twice continuously-differentiable functions $\phi(x)$ with
restricted support in $(a, b)$ such that $\phi(x) \ge 0$ for all $x
\in (a, b)$.  If $f(x)$ is also twice continuously-differentiable,
then $f'' \ge 0$ on $(a, b)$ in the weak sense if and only if 
$f''(x) \ge 0$ for all $x \in (a, b)$, basically by integrating
by parts twice.  For any continuous real-valued function $f(x)$
on $(a, b)$, one can show that $f(x)$ is convex if and only if 
$f'' \ge 0$ on $(a, b)$ in the weak sense.

\subsection{Subharmonic functions, 2}
\label{subsection on subharmonic functions, 2}

	Let $U$ be a nonempty open subset of ${\bf R}^n$, and let
$f(x)$ be a real-valued continuous function on $U$.  We say that
$f(x)$ satisifies the weak sub-mean value property on $U$ if for
each $x \in U$ there is a positive real number $r(x)$ such that
\begin{equation}
	{\average}_{\overline{B}(x,r)} f \ge f(x)
\end{equation}
for all $r \in (0, r(x))$.  To be more precise, this is a sub-mean
value property for averages over balls, and we shall also consider
versions with averages over spheres.

	By definition, if a continuous real-valued function satisfies
the weak sub-mean value property on an open set, it also satisfies
this property when restricted to any open subset of this set.  If $U$
is a nonempty connected open subset of ${\bf R}^n$, $f(x)$ is a
real-valued continuous function on $U$ which satisfies the weak
sub-mean value property, and $p$ is an element of $U$ at which $f(x)$
attains its maximum, so that
\begin{equation}
	f(x) \le f(p),
\end{equation}
then $f(x)$ is constant on $U$.  As usual, this is because $\{x \in U
: f(x) = f(p)\}$ is then an open subset of $U$, by the weak sub-mean
value property, while $\{x \in U : f(x) < f(p)\}$ is an open subset of
$U$ by continuity.

	Let $U$ be a nonempty open subset of ${\bf R}^n$, let $f(x)$
be a real-valued continuous function on $U$, let $z$ be an element of
$U$, and let $r$ be a positive real number such that
\begin{equation}
	r < \dist(z, {\bf R}^n \backslash U).
\end{equation}
Let $h(x)$ be the continuous real-valued function on $\overline{B}(z,
r)$ which is equal to $f(x)$ when $|x - z| = r$ and which is harmonic
on $B(z, r)$.  Thus $f(x) - h(x)$ satisfies the weak sub-mean value
property on $B(z, r)$, and since $f(x) - h(x)$ is equal to $0$ on the
boundary, it follows that $f(x) - h(x) \le 0$ on $B(z, r)$, which is
to say that
\begin{equation}
	f(x) \le h(x)
\end{equation}
when $|x - z| < r$.

	Because $h(x)$ is harmonic on $B(z, r)$, we have that
\begin{equation}
	\int_{B(z, r)} h(x) \, b(x) \, dx = h(z)
\end{equation}
whenever $b(x)$ is a continuous real-valued function with restricted
support in $B(z, r)$ which is radial around $z$, in the sense that it
depends only on $|x - z|$, and which satisfies
\begin{equation}
	\int_{B(z, r)} b(x) \, dx = 1.
\end{equation}
As a result,
\begin{equation}
	{\average}_{\partial B(z, t)} h = h(z)
\end{equation}
when $0 < t < r$, where $\partial B(z, t) = \{x \in {\bf R}^n : |x -
z| = t\}$.  This equation also holds when $t = r$, since $h(x)$ is
continuous on $\overline{B}(z, r)$.

	It follows that
\begin{equation}
	{\average}_{\partial B(z, r)} f 
		\ge {\average}_{\partial B(z, t)} f
\end{equation}
when $0 < t < r$, and that
\begin{equation}
	{\average}_{\partial B(z, r)} f \ge f(z).
\end{equation}
This works for all radii $r < \dist(z, {\bf R}^n \backslash U)$, and
by integrating over $r$ one gets that if $\beta(x)$ is a nonnegative
real-valued continuous function with restricted support in $U$ which
is radial about $z$, so that $\beta(x)$ only depends on $|x - z|$, and
which satisfies
\begin{equation}
	\int_U \beta(x) \, dx = 1,
\end{equation}
then
\begin{equation}
	\int_U f(x) \, \beta(x) \, dx \ge f(z).
\end{equation}
Also,
\begin{equation}
	{\average}_{\overline{B}(z, r)} f \ge f(z)
\end{equation}
when $z \in U$ and $0 < r < \dist(z, {\bf R}^n \backslash U)$.

\subsection{Subharmonic functions, 3}
\label{subsection on subharmonic functions, 3}

	Let $U$ be a nonempty open subset of ${\bf R}^n$.  A
real-valued continuous function $f(x)$ is said to be
\emph{subharmonic}\index{subharmonic functions on ${\bf R}^n$} if it
satisfies the weak sub-mean value property on $U$, which is then
equivalent to seemingly-stronger conditions, as in the previous
subsection.  If $f(x)$ is twice continuously-differentiable on $U$,
then $\Delta f(x) \ge 0$ for all $x \in U$ in this case.

	Conversely, suppose that $f(x)$ is a real-valued function on a
nonempty open subset $U$ of ${\bf R}^n$ such that $f$ is twice
continuously-differentiable and $\Delta f(x) \ge 0$.  Actually, we can
weaken this a bit by assuming that $f(x)$ is a real-valued continuous
function on $U$ such that
\begin{equation}
	\liminf_{r \to 0} \, r^{-2} \,
	 \biggl({\average}_{\overline{B}(x, r)} f - f(x)\biggr) 
		\ge 0
\end{equation}
for all $x \in U$, or even
\begin{equation}
	\limsup_{r \to 0} \, r^{-2} \,
	 \biggl({\average}_{\overline{B}(x, r)} f - f(x)\biggr) 
		\ge 0
\end{equation}
for all $x \in U$.  In this case one can check that
\begin{equation}
	f_\epsilon(x) = f(x) + \epsilon \, |x|^2
\end{equation}
satisfies the weak sub-mean value property on $U$ for all $\epsilon >
0$.  Hence
\begin{equation}
	{\average}_{\overline{B}(x, r)} f_\epsilon \ge f_\epsilon(x)
\end{equation}
when $x \in U$ and $0 < r < \dist(x, {\bf R}^n \backslash U)$, and
therefore
\begin{equation}
	{\average}_{\overline{B}(x, r)} f \ge f(x)
\end{equation}
for the same $x$, $r$, so that $f(x)$ is a subharmonic function on
$U$.

	A real-valued continuous function $f(x)$ on a nonempty open
subset $U$ of ${\bf R}^n$ is said to be \emph{subharmonic in the weak
sense}\index{weakly subharmonic functions on an open subset of ${\bf
R}^n$} if
\begin{equation}
	\int_U f(x) \, \Delta \phi(x) \, dx \ge 0
\end{equation}
for all nonnegative real-valued twice continuously-differentiable
functions $\phi(x)$ with restricted support in $U$.  If $f(x)$ is
twice continuously-differentiable, then it is easy to check that this
is equivalent to $\Delta f(x) \ge 0$ for all $x \in U$, by integration
by parts.  It is not too difficult to extend this argument to show
that a subharmonic function is subharmonic in the weak sense.

	Conversely, a real-valued continuous function $f(x)$ on $U$
which is subharmonic in the weak sense is subharmonic.  This is
analogous to the fact that a function which is harmonic in the weak
sense satisfies the mean value property.  As before, one can express
differences of local radial averages of a function $f$ as integrals of
$f$ times $\Delta \phi$ for twice continuously-differentiable
functions with restricted support in $U$, and now one should be a bit
more careful to have $\phi \ge 0$ too.

	Of course a real-valued continuous function $f(x)$ on an open
subset $U$ of ${\bf R}^n$ is harmonic if and only if $f(x)$ and
$-f(x)$ are subharmonic.  If $f_1$, $f_2$ are subharmonic functions on
$U$ and $a_1$, $a_2$ are nonnegative real numbers, then
\begin{equation}
	a_1 \, f_1 + a_2 \, f_2
\end{equation}
is also a subharmonic function on $U$.  Moreover,
\begin{equation}
	\max(f_1, f_2)
\end{equation}
is a subharmonic function on $U$ in this case.

\subsection{Holomorphic functions and subharmonicity}
\label{subsection on holomorphic functions and subharmonicity}

	Let $U$ be a nonempty open subset of ${\bf R}^n$, and
let $f(x)$ be a real-valued continuous function on $U$.  If
$\alpha(t)$ is a real-valued convex function on ${\bf R}$, then
\begin{equation}
	\alpha\biggl({\average}_{\overline{B}(x, r)} f\biggr)
		\le {\average}_{\overline{B}(x, r)} \alpha \circ f
\end{equation}
for all $x \in U$ and $0 < r < \dist(x, {\bf R}^n \backslash U)$.
In other words, the value of a convex function at an average of
some numbers is less than or equal to the average of the values
of the convex function at the same numbers.

	As a consequence, if $f(x)$ is a real-valued harmonic function
on $U$ and $\alpha(t)$ is a convex function on the real line, then the
composition $\alpha(f(x))$ is a real-valued continuous function on $U$
which is subharmonic.  Moreover, if $f(x)$ is a real-valued continuous
subharmonic function on $U$ and $\alpha(t)$ is a monotone increasing
convex function on the real line, then the composition $\alpha(f(x))$
is a real-valued continuous subharmonic function on $U$.  For
instance, in the first case one can take $\alpha(t) = |t|^p$ and in
the second case one can take $\alpha(t) = 0$ when $t < 0$, $\alpha(t)
= t^p$ when $t \ge 0$, where $p$ is a real number such that $p \ge 1$.

	In fact, if $U$ is a nonempty open subset of ${\bf C}$, and
$f(z)$ is a complex-valued holomorphic function on $U$, then
$|f(z)|^p$ is subharmonic for all $p > 0$.  One also has that $\log
|f(z)|$ is subharmonic, if one extends the notion of subharmonicity to
allow for functions which take the value $-\infty$, or alternatively
the real-valued function $\max(\log|f(z)|, \lambda)$ is subharmonic
for every real number $\lambda$.  Indeed, these assertions can be
verified using the fact that $\log |f(z)|$ is harmonic on the set
where $f(z) \ne 0$.

	Now suppose that $U$ is a nonempty open subset of ${\bf R}^n$,
and that $f(x)$ is a $\mathcal{C}(n)$-valued function on $U$.  We can
define $|f(x)|$ by identifying $\mathcal{C}(n)$ with a Euclidean
space, using the standard basis for the Clifford algebra.  It turns
out that if $f(x)$ is Clifford holomorphic, then $|f(x)|^p$ is
subharmonic for a range of $p$'s which includes some $p$'s strictly
less than $1$.  There are general results of this type for first-order
constant coefficient systems of partial differential equations which
are ``elliptic'' in the sense that their homogeneous solutions are
harmonic, and there are various situations where this applies with
various ranges of $p$.  See \cite{Stein, SW}.

\section{Miscellaneous}
\label{section on miscellaneous matters}
\setcounter{equation}{0}

\subsection{Linear fractional transformations}
\label{subsection on linear fractional transformations}

	A \emph{polynomial}\index{polynomials on ${\bf C}$}
on ${\bf C}$ is a function $P(z)$ of the form
\begin{equation}
	a_n \, z^n + a_{n-1} \, z^{n-1} + \cdots + a_0,
\end{equation}
where $n$ is a nonnegative integer and $a_0, \ldots, a_n$ are complex
numbers.  If $n \ge 1$ and $a_n \ne 0$, or if $n = 0$, then we say
that $P(z)$ has degree $n$.  If $P(z)$ is a polynomial of degree $n
\ge 1$, then every complex number is attained as a value of $P(z)$
exactly $n$ times if one counts multiplicities in an appropriate
manner.

	Let us write $\widehat{\bf C}$ for the \emph{Riemann sphere},
which consists of the complex numbers ${\bf C}$ together with an extra
``point at infinity'' denoted $\infty$.  This is topologically
equivalent to the usual $2$-dimensional sphere in the obvious manner.
We also make the conventions that $1/\infty = 0$, $1/0 = \infty$,
$a + \infty = \infty$ when $a \in {\bf C}$, $a \cdot \infty = \infty$
when $a \in {\bf C}$, $a \ne 0$, etc.

	If $P(z)$ is a polynomial on ${\bf C}$, then $P$ can be 
extended to a continuous mapping from the Riemann sphere in a 
natural way, namely, $P(\infty) = \infty$ when $P(z)$ has positive
degree, and $P$ takes the same value at $\infty$ as at other points
when $P$ is constant.  More generally we can consider \emph{rational
functions}\index{rational functions of a complex behavior} of the
form
\begin{equation}
	R(z) = \frac{P(z)}{Q(z)},
\end{equation}
where $Q(z)$ is a polynomial which is not identically equal to $0$.
To be more precise, $R(z)$ is then defined as a complex number when
$z$ is a complex number and $Q(z) \ne 0$, while if $Q(z) = 0$ or $z =
\infty$ one can follow the usual conventions, so that the rational
function remains continuous as a mapping from the Riemann sphere to
itself.

	Assuming that the polynomials $P(z)$, $Q(z)$ have no common
factors, then the degree of the corresponding rational function $R(z)$
is defined to be the maximum of the degrees of $P(z)$ and $Q(z)$.  If
$R(z)$ is a rational function of degree $n \ge 1$, then each element
of the Riemann sphere occurs as a value of $R(z)$ $n$ times, with
appropriate multiplicities.  This is not too difficult to see.

	A \emph{linear fractional transformation}\index{linear
fractional transformation} on the Riemann sphere is a rational
function of degree $1$.  Explicitly, this is a mapping which can be
expressed as
\begin{equation}
	\frac{a \, z + b}{c \, z + d},
\end{equation}
where $a$, $b$, $c$, $d$ are complex numbers such that
\begin{equation}
	a \, d - b \, c \ne 0.
\end{equation}
This condition ensures that the denominator in the linear fractional
transformation is not identically equal to $0$, and that the numerator
and denominator are not simply constant multiples of each other.
Notice that if we multiply $a$, $b$, $c$, and $d$ by the same nonzero
complex number, then this admissibility condition is still satisfied, 
and that the corresponding linear fractional transformation on the
Riemann sphere is in fact the same as the original one.

	In fact, we can think of the complex numbers $a$, $b$, $c$,
$d$ as entries of a $2 \times 2$ matrix
\begin{equation}
	\biggl({a \ b \atop c \ d}\biggr),
\end{equation}
whose determinant is then equal to $a \, c - b \, d$, so that the
admissibility condition is equivalent to the requirement that this
matrix be invertible.  One can check that the composition of two
linear fractional transformations is again a linear fractional
transformations, and that this corresponds to multiplying the
associated matrices in the usual way.  In particular, the inverse
of a linear fractional transformation can be given by the inverse
of the associated matrix.

	It is helpful to consider some special cases of linear
fractional transformations, of the form
\begin{equation}
	a \, z
\end{equation}
for a nonzero complex number $a$,
\begin{equation}
	z + b
\end{equation}
for any complex number $b$, and
\begin{equation}
	\frac{1}{z}.
\end{equation}
One can check that every linear fractional transformation can be
expressed as a combination of these basic types.  Of course this
can also be viewed in terms of $2 \times 2$ matrices.

	Let us consider this in another way, which extends to higher
dimensions.  Namely, if $n$ is a positive integer, let ${\bf CP}^n$
denote the space of complex lines through the origin in ${\bf
C}^{n+1}$, the space of $(n+1)$-tuples of complex numbers.
In other words, we can start with ${\bf C}^{n+1} \backslash \{0\}$,
and identify two elements $z$, $w$ of this space when there is a
nonzero complex number $\lambda$ such that $w = \lambda \, z$.

	If $A$ is an invertible complex linear mapping on ${\bf
C}^{n+1}$, then $A$ takes the origin $0$ in ${\bf C}^{n+1}$ to itself,
and it takes complex lines in ${\bf C}^{n+1}$ to themselves, and thus
induces a transformation on ${\bf CP}^n$.  Multiplication of $A$ by a
nonzero complex number does not change the induced transformation on
${\bf CP}^n$.  Compositions of invertible linear mappings on ${\bf
C}^{n+1}$ correspond to compositions of the induced transformations on
${\bf CP}^n$, and in particular the inverse of an induced
transformation on ${\bf CP}^n$ is induced by the inverse of the
original linear mapping on ${\bf C}^{n+1}$.

	There is a natural way to embedd ${\bf C}^n$ into ${\bf
CP}^n$, which is to take an element $z$ of ${\bf C}^n$ and extend it
to a nonzero element of ${\bf C}^{n+1}$ by adding an $(n+1)$th
coordinate which is set equal to $1$.  In other words, we identify
${\bf C}^n$ with the affine hyperspace in ${\bf C}^{n+1}$ of points
whose last coordinate is equal to $1$, and then each point in this
affine hyperspace determines a complex line in ${\bf C}^{n+1}$ through
the origin, namely the complex line that passes through that point.
This accounts for all the complex lines in ${\bf C}^{n+1}$ except for
those which lie in the linear subspace of points whose last coordinate
is equal to $0$.

	Roughly speaking then, ${\bf C}^n$ looks like ${\bf C}^n$
together with a copy of ${\bf CP}^{n-1}$, at least when $n \ge 2$.
When $n = 1$, the embedding of ${\bf C}$ into ${\bf CP}^1$ accounts
for all the lines in ${\bf C}^2$ through the origin except for one
line, which is the set of points with second coordinate equal to $0$.
Thus ${\bf CP}^1$ looks like ${\bf C}$ with one additional element,
which gives another way to think about the Riemann sphere.

	Linear fractional transformations on the Riemann sphere are
then the same as transformations on ${\bf CP}^1$ induced by invertible
complex linear transformations on ${\bf C}^2$.  To be more explicit,
an invertible linear transformation on ${\bf C}^2$ can be given as
\begin{equation}
	(z, w) \mapsto (a \, z + b \, w, c \, z + d \, w),
\end{equation}
where again $a$, $b$, $c$, $d$ are complex numbers such that 
$a \, c - b \, d \ne 0$.  This leads to the earlier formula through
the identifications described above.

	There is another and rather different way to extend the idea
of linear fractional transformations to higher dimensions.  For this
we think of ${\bf R}^2$ as being identified with ${\bf C}$ in the
usual way.  Now we work on ${\bf R}^n$, and on $\widehat{{\bf R}^n}
\simeq {\bf R}^n \cup \{\infty\}$, which can be identified
topologically with an $n$-dimensional sphere.

	To describe this generalization, one can go back to the idea
of building blocks.  Namely, one can use translations, dilations,
orthogonal transformations, and inversion about the unit sphere in
${\bf R}^n$ as the basic building blocks.  One then gets a nice family
of transformations on $\widehat{{\bf R}^n}$ generated by these.

	Actually, this would correspond to also allowing complex
conjugation as a transformation on the Riemann sphere.  Alternatively,
in $n$ dimensions, one can add the requirement that the
transformations preserve orientations.  At any rate, this
generalization is closely connected to conformal geometry in $n$
dimensions.

\subsection{Power series and the exponential function}
\label{subsection on power series and the exponential function}

	Consider a power series
\begin{equation}
	\sum_{n=0}^\infty a_n \, z^n,
\end{equation}
where the $a_n$'s are complex numbers and $z$ is a complex variable.
We shall be interested in the set of $z \in {\bf C}$ for which this
series converges.  Of course this series converges trivially when $z =
0$.

	If the series converges for some $z_0 \in {\bf C}$, then we
have that the sequence $\{a_n \, z_0^n\}_{n=0}^\infty$ tends to $0$,
and is bounded in particular, so that there is a nonnegative real
number $C$ such that
\begin{equation}
	|a_n| \, |z_0|^n \le C
\end{equation}
for all $n$.  Assume that $z_0 \ne 0$ and that $z$ is a complex number
such that
\begin{equation}
	|z| < |z_0|.
\end{equation}
Under these conditions, we have that
\begin{equation}
	|a_n| \, |z|^n \le C \, \biggl(\frac{|z|}{|z_0|}\biggr)^n,
\end{equation}
and hence $\sum_{n=0}^\infty a_n \, z^n$ converges absolutely by 
comparison with a convergent geometric series.

	Using this simple remark, one can check that either the power
series converges only for $z = 0$, or there is a positive real number
$R$ such that the series converges absolutely when $|z| < R$ and does
not converge when $|z| > R$, or the series converges absolutely for
all $z \in {\bf C}$.  In the first case we can set $R = 0$, and in the
third case we can set $R = \infty$.  We call $R$ the radius of
convergence of the power series.

	If the power series converges absolutely for some complex
number $z_0$, then it also converges absolutely for all complex
numbers $z$ such that $|z| \le |z_0|$.  In general, if a power series
has a radius of convergence $R$ which is positive and finite, then
the series might converge for no complex numbers $z$ with $|z| = R$,
or for some of them, or for all of them.  For instance, the series
\begin{equation}
	\sum_{n=0}^\infty z^n
\end{equation}
has radius of convergence equal to $1$ and converges for no $z \in
{\bf C}$ with $|z| = 1$,
\begin{equation}
	\sum_{n=1}^\infty \frac{z^n}{n}
\end{equation}
has radius of convergence $1$ and converges for $z \in {\bf C}$ with
$|z| = 1$ and $z \ne 1$, and does not converge when $z = 1$, and
\begin{equation}
	\sum_{n=1}^\infty \frac{z^n}{n^2}
\end{equation}
has radius of convergence $1$ and converges absolutely for all $z \in
{\bf C}$ with $|z| \le 1$.

	Let
\begin{equation}
	\sum_{n=0}^\infty A_n, \quad \sum_{n=0}^\infty B_n
\end{equation}
be two series of complex numbers, and define an associated series
\begin{equation}
	\sum_{n=0}^\infty C_n
\end{equation}
by
\begin{equation}
	C_n = \sum_{j=0}^n A_j \, B_{n-j}.
\end{equation}
This is called the Cauchy product of the original two series.
Formally, the sum of the $C_n$'s is equal to the product of the sums
of the $A_n$'s and the $B_n$'s, and indeed this is the case if the
$A_n$'s and $B_n$'s are equal to $0$ for all but finitely many $n$.
If the original series are power series, with $A_n = a_n \, z_n$ and
$B_n = b_n \, z^n$, then $C_n = c_n \, z^n$, with $c_n = \sum_{j=0}^n
a_j \, b_{n-j}$.

	It is not too difficult to show that if $\sum A_n$, $\sum B_n$
converge absolutely, then the Cauchy product series $\sum C_n$ also
converges absolutely, and the sum of the Cauchy product is equal to
the products of the sums of the initial series.  One can also show
that if both series converge and at least one converges absolutely,
then the Cauchy product converges, and the sum of the Cauchy product
series is equal to the product of the sums of the original series.
It is not true in general that the Cauchy product series converges
when the two initial series converge.

	There is a nice result though which says that if the two
initial series converge and the Cauchy product series converges, then
the sum of the Cauchy product series is equal to the product of the
sums of the two initial series.  To show this, one can use the notion
of \emph{Abel summability}.  Basically, it is convenient to view the
series as power series.

	Specifically, for a given series $\sum_{n=0}^\infty A_n$,
the associated Abel sums are defined by
\begin{equation}
	\sum_{n=0}^\infty r^n \, A_n,
\end{equation}
where $r$ is a positive real number such that $r < 1$.  We assume that
these series converge for all $r \in (0, 1)$, which is the same as
saying that they converge absolutely for all $r \in (0,1)$, as before.
This holds under modest assumptions on the size of the $A_n$'s.

	The series $\sum_{n=0}^\infty A_n$ is said to be Abel summable
if the Abel sums exist for all $r \in (0,1)$, and if the limit of the
Abel sums exist as $r \to 1-$.  This limit is then called the Abel sum
of the series $\sum_{n=0}^\infty A_n$.  For instance, if the $A_n$'s
are nonnegative real numbers, then the Abel sums are monotone
increasing as a function of $r$, and it is easy to see that the series
is Abel summable if and only if it converges in the usual sense,
and the Abel sum is equal to the ordinary sum.

	If $\sum A_n$ converges absolutely, then it is again easy to
see that $\sum A_n$ is Abel summable, and that the Abel sum is equal
to the ordinary sum.  It turns out that this works in general when
$\sum A_n$ converges.  One can approach this by expressing the Abel
sums in terms of the partial sums of the original series, noting also
that the Abel sums converge for all $r \in (0,1)$ because the $A_n$'s
are bounded.

	There are also series which are Abel summable and which do not
converge themselves.  A basic family of examples is given by the
series
\begin{equation}
	\sum_{n=0}^\infty z^n,
\end{equation}
where $z$ is a complex number such that $|z| = 1$ and $z \ne 1$.
When $r \in (0, 1)$, we have that
\begin{equation}
	\sum_{n=0}^\infty r^n \, z^n = \frac{1}{1 - r \, z},
\end{equation}
and this converges to $(1 - z)^{-1}$ as $r \to 1-$.

	Now suppose that $\sum A_n$, $\sum B_n$ are two convergent
series of complex numbers.  In particular their terms are bounded, and
the Cauchy product series $\sum C_n$ has the property that $C_n =
O(n+1)$, say, which is sufficient to ensure that $\sum r^n \, C_n$
converges for all $r \in (0,1)$.  If we assume that $\sum C_n$
converges, then the sum of this series is the same as the limit of the
Abel sums $\sum r^n \, C_n$, and these Abel sums are equal to the
product of the Abel sums $\sum r^n \, A_n$, $\sum r^n \, B_n$, since
we have absolute convergence when $r \in (0,1)$, and it follows
that $\sum C_n$ is equal to $\sum A_n$ times $\sum B_n$.

	Suppose that we have a power series
\begin{equation}
	\sum_{n=0}^\infty \alpha_n \, z^n
\end{equation}
which converges for $z \in {\bf C}$ with $|z| < R$ for some $R > 0$.
A basic result states that if $0 < r < R$, then the partial sums of
this series converge uniformly to the whole sum for $z \in {\bf C}$
such that $|z| \le r$.  In particular, it follows that the power
series defines a continuous function on the disk 
\begin{equation}
	\{z \in {\bf C} : |z| < R\}.
\end{equation}
Under these conditions we also have that the series
\begin{equation}
	\sum_{n=0}^\infty n \, \alpha_n \, z^n
\end{equation}
converges for $|z| < R$.  The original series defines a holomorphic
function $f(z)$ on the open disk in the complex plane with center $0$
and radius $R$, and this series represents the complex derivative
$f'(z)$ of $f(z)$.

	Now let us look at a specific power series.  Namely,
consider the series
\begin{equation}
	\sum_{n=0}^\infty \frac{z^n}{n!},
\end{equation}
where as usual $n!$ denotes ``$n$ factorial'', the product of the
integer from $1$ to $n$, which is interpreted as being equal to $1$
when $n = 0$.  For each positive integer $k$, we have that
\begin{equation}
	\frac{1}{n!} \le \frac{1}{(k-1)!} \frac{1}{k^{n-k}}
\end{equation}
when $n \ge k$, and one can use this to show that our power series
converges for all $z$ in ${\bf C}$.

	Of course this is the series for the exponential function
$\exp (z)$.  With this definition one can see directly that
\begin{equation}
	\exp (z + w) = \exp(z) \, exp(w)
\end{equation}
for all $z, w \in {\bf C}$.  This is an example of the general 
result for Cauchy products of absolutely convergent series,
as one can check using the binomial theorem.

	Clearly $\exp (0) = 1$.  When $z$ is a positive real number,
$\exp (z)$ is a positive real number greater than $1$.  If $z$ is a
negative real number, $\exp (z) = 1/\exp(-z)$ is a positive real
number less than $1$.

	For any complex number $z$, 
\begin{equation}
	\overline{\exp (z)} = \exp (\overline{z}),
\end{equation}
where $\overline{a}$ denotes the complex conjugate of a complex number
$a$.  In particular, if $z = x + i \, y$, where $x$, $y$ are real
numbers, then
\begin{equation}
	|\exp (z)| = \exp (x),
\end{equation}
which is to say that
\begin{equation}
	|\exp (i \, y)| = 1.
\end{equation}
In fact,
\begin{equation}
	\exp (i \, y) = \cos y + i \, \sin y,
\end{equation}
as one can see from the power series expansions.

	The exponential function is a holomorphic function of $z$
whose complex derivative is equal to itself, just as the ordinary
derivative of $\exp (x)$ for $x \in {\bf R}$ is equal to itself, as
one can see from the power series.  Note that one can also derive
relations for the derivatives of the real and imaginary parts of $\exp
(i \, y)$ as a function on the real line which are compatible with the
formula in terms of sines and cosines above.  The fact that the
exponential of a sum is equal to the product of the individual
exponentials also contains standard addition formulas for sines and
cosines.

	As a consequence of the trigonometric relations, we 
get that
\begin{equation}
	\exp (2 \pi \, i \, n) = 1
\end{equation}
when $n$ is an integer, and that these are the only times when $\exp
(z) = 1$.  Also, $\exp (z) \ne 0$ for all $z \in {\bf C}$, and if
$\zeta$ is a nonzero complex number, then there are $z$'s in ${\bf C}$
such that
\begin{equation}
	\exp (z) = \zeta.
\end{equation}
Specifically, the real part $x$ of $z$ is uniquely determined by the
condition that $\exp (x) = |\zeta|$, while the imaginary part $y$ of
$z$ should be taken to be an ``angle'' for $\zeta$, as in polar
coordinates for ${\bf R}^2 \simeq {\bf C}$, which is determined up to
adding integer multiples of $2 \pi$.

	If $U$ is a nonempty open subset of ${\bf C} \backslash \{0\}$
and $\lambda(\zeta)$ is a continuous complex-valued function on $U$,
then we say that $\lambda(\zeta)$ is a \emph{branch of the logarithm}
on $U$ if
\begin{equation}
	\exp (\lambda(\zeta)) = \zeta
\end{equation}
for all $\zeta \in U$.  Assuming that $\lambda_1(\zeta)$,
$\lambda_2(\zeta)$ are branches of the logarithm on the same nonempty
open subset $U$ of ${\bf C} \backslash \{0\}$, we get that
\begin{equation}
	\lambda_1(\zeta) - \lambda_2(\zeta) 
			\in 2 \pi \, i \, {\bf Z}
\end{equation}
for all $\zeta \in U$.  Under the additional hypothesis that $U$ is
connected, it follows that there is a single integer $n$ such that
\begin{equation}
	\lambda_1(\zeta) - \lambda_2(\zeta) = 2 \pi \, i \, n
\end{equation}

	As above, if $\lambda(\zeta)$ is a branch of the logarithm
on a nonempty open subset $U$ of ${\bf C} \backslash \{0\}$, then
\begin{equation}
	\re \lambda(\zeta) = \log |\zeta|,
\end{equation}
where the logarithm $\log \rho$ of a positive real number $\rho$
is the real number $r$ characterized by $\exp (r) = \rho$.
Also, $\lambda(\zeta)$ is a holomorphic function on $U$.  The
complex derivative of $\lambda(\zeta)$ is given by
\begin{equation}
	\lambda'(\zeta) = \frac{1}{\zeta}.
\end{equation}

	The function $1/\zeta$, up to a nonzero constant multiple, is
a fundamental solution for $\partial / \partial \overline{z}$ on ${\bf
C}$.  On ${\bf R}^n$ for general $n$ one has the Clifford-valued
function
\begin{equation}
	\frac{\sum_{j=1}^n x_j \, e_j}{|x|^n},
\end{equation}
a nonzero constant multiple of which is a fundamental solution for
$\mathcal{D}_L$, $\mathcal{D}_R$.  Integrals of this function are also
related to measurements of angular increments, as in the complex
plane.

\subsection{Linear algebra in ${\bf R}^n$, continued}
\label{subsection on linear algebra in R^n, continued}

	Let $U$ be a nonempty open subset of ${\bf R}^n$, and let
$f(x)$ be a real-valued twice continuously-differentiable function on
$U$.  For each $x \in U$, we get the matrix of second partial
derivatives of $f$,
\begin{equation}
	\frac{\partial^2}{\partial x_j \partial x_l}, 
					\quad 1 \le j, l \le n,
\end{equation}
and a theorem of advanced calculus states that this matrix is
symmetric.  The Laplacian $\Delta f(x)$ of $f$ at $x$ is the same as
the trace of this matrix.

	If $T$ is a linear transformation on ${\bf R}^n$ and $v$ is a
nonzero vector in ${\bf R}^n$, then $v$ is said to be an
\emph{eigenvector}\index{eigenvector of a linear transformation on
${\bf R}^n$} for $T$ with \emph{eigenvalue}\index{eigenvalues of a
linear transformation on ${\bf R}^n$} $\lambda \in {\bf R}$ if
\begin{equation}
	T(v) = \lambda \, v.
\end{equation}
If $T$ is a symmetric linear transformation on ${\bf R}^n$, 
and if $v$ is an eigenvector of $T$ with eigenvalue $\lambda$,
then $T$ maps every vector $w$ in ${\bf R}^n$ orthogonal to $v$
to a vector orthogonal to $v$.  Indeed,
\begin{eqnarray}
	\langle T(w), v \rangle & = & \langle w, T(v) \rangle	\\
		& = & \lambda \, \langle w, v \rangle = 0.
						\nonumber
\end{eqnarray}

	In general, a linear transformation $T$ on ${\bf R}^n$ is said
to be \emph{diagonalizable}\index{diagonalizable linear
transformations on ${\bf R}^n$} if there is a basis of ${\bf R}^n$
consisting of eigenvectors of $T$, which is equivalent to saying that
there is an invertible linear transformation $A$ on ${\bf R}^n$ so
that the matrix associated to $A \, T \, A^{-1}$ is diagonal.  If $T$
is a symmetric linear transformation on ${\bf R}^n$, then there is an
orthonormal basis of ${\bf R}^n$ consisting of eigenvectors of $T$,
which is equivalent to saying that there is an orthogonal linear
transformation $A$ on ${\bf R}^n$ such that the matrix associated to
$A \, T \, A^{-1}$ is diagonal.  To find an eigenvector $v$ for $T$,
one can maximize
\begin{equation}
	\langle T(v), v \rangle
\end{equation}
among vectors $v$ such that $|v| = 1$, and then one can repeat the
process to get an orthonormal basis of eigenvectors.

	Suppose that $T$ is any linear transformation on ${\bf R}^n$.
Then $T^* \, T$ is a nonnegative symmetric linear operator on
${\bf R}^n$, and there is a nonnegative symmetric linear operator
$\widehat{T}$ on ${\bf R}^n$ such that
\begin{equation}
	\widehat{T}^2 = T.
\end{equation}
Namely, there is an orthonormal basis of ${\bf R}^n$ consisting of
eigenvectors of $T^* \, T$, with nonnegative eigenvalues, and one can
define $\widehat{T}$ so that the same basis is a basis of eigenvectors
and the eigenvalues are nonnegative square roots of those of $T$.

	Observe that 
\begin{equation}
	|\widehat{T}(v)| = |T(v)|
\end{equation}
for all $v \in {\bf R}^n$, since
\begin{eqnarray}
	|\widehat{T}(v)|^2 & = & 
		\langle \widehat{T}(v), \widehat{T}(v) \rangle  \\
		& = & \langle \widehat{T}^2(v), v \rangle
						\nonumber	\\
		& = & \langle (T^* \, T)(v), v \rangle
		  = \langle T(v), T(v) \rangle = |T(v)|^2.
						\nonumber
\end{eqnarray}
Using this one can show that there is an orthogonal linear
transformation $A$ on ${\bf R}^n$ such that
\begin{equation}
	T = A \circ \widehat{T}.
\end{equation}
This is called a \emph{polar decomposition}\index{polar decomposition
of a linear transformation on ${\bf R}^n$} of $T$, although one might
prefer a decomposition of $T$ in the other order, as the composition
of a nonnegative symmetric operator with an orthogonal transformation,
and this can be obtained by expressing $T^*$ as the composition of an
orthogonal transformation and a nonnegative symmetric operator, and
taking the adjoint afterwards.  For this polar decomposition, one can
simply take $A = T \, \widehat{T}^{-1}$ when $T$ is invertible, so
that $\widehat{T}$ is invertible, and otherwise one has some freedom
in choosing $A$ on the orthogonal complement of the image of
$\widehat{T}$.

	Now let us look some more at orthogonal transformations in
general.  If $v_1, \ldots, v_n$ and $w_1, \ldots, w_n$ are orthogonal
bases on ${\bf R}^n$, then there is a unique linear transformation $A$
on ${\bf R}^n$ such that $A(v_j) = w_j$ for $1 \le j \le n$, and $A$
is an orthogonal transformation.  Conversely, if $A$ is an orthogonal
linear transformation on ${\bf R}^n$ and $v_1, \ldots, v_n$ is an
orthonormal basis of ${\bf R}^n$, then $A(v_1), \ldots, A(v_n)$ is
also an orthonormal basis of ${\bf R}^n$.

	If $A$ is an orthogonal linear transformation on ${\bf R}^n$
and $v_1, \ldots, v_n$ is an orthonormal basis of ${\bf R}^n$, then
$A$ is almost determined by its values on $v_1, \ldots, v_{n-1}$.
Namely, $A(v_n)$ should be a unit vector orthogonal to $A(v_1),
\ldots, A(v_{n-1})$, and there are only two choices for this vector,
which are negatives of each other.  If $A$ is a rotation, which is to
say an orthogonal linear transformation with determinant $1$, so that
$A$ preserves orientations on ${\bf R}^n$, then $A$ is determined by
its values on $v_1, \ldots, v_{n-1}$.

	Suppose that $u_1, \ldots, u_k$ are orthonormal vectors in
${\bf R}^n$ with $k \le n - 2$, or simply that $n \ge 2$ if no such
vectors are specified.  If $v$, $w$ are unit vectors in ${\bf R}^n$
which are orthogonal to $u_1, \ldots, u_k$, then one can make a
continuous deformation of $v$ to $w$ through unit vectors in ${\bf
R}^n$ which are orthogonal to $u_1, \ldots, u_k$.  Indeed, there is a
continuous deformation of the identity operator on ${\bf R}^n$
through orthogonal linear transformations on ${\bf R}^n$ which
fixes $u_j$ for $j = 1, \ldots, k$ and which deforms $v$ to $w$.

	Using this, one can show that any two rotations on ${\bf R}^n$
can be continuously deformed to each other through orthogonal linear
transformations.  For that matter, any two invertible linear
transformations on ${\bf R}^n$ whose determinants have the same sign
can be continuously deformed to each other.  If $A_1$, $A_2$ are two
invertible linear transformations on ${\bf R}^n$ whose determinants
have opposite signs, then any continuous deformation from $A_1$ to
$A_2$ through linear transformations on ${\bf R}^n$ must pass through
at least one linear transformation with determinant.

	Recall that a function $f(x)$ on ${\bf R}^n$ is said to be
radial if it depends only on $|x|$.  This is equivalent to saying that
\begin{equation}
	f(A(x)) = f(x)
\end{equation}
for all orthogonal linear transformations $A$ on ${\bf R}^n$.  For if
$x$, $y$ are elements of ${\bf R}^n$ such that $|x| = |y|$, then there
is an orthogonal linear transformation $A$ on ${\bf R}^n$ such that
$A(x) = y$.

	Suppose that $U$ is a nonempty open subset of ${\bf R}^n$ and
that $f$ is a twice continuously-differentiable function on $U$.  If
$A$ is an orthogonal linear transformation on ${\bf R}^n$, then $f
\circ A$ is a twice continuously-differentiable function on
$A^{-1}(U)$.  One can check that the Laplacian of $f \circ A$ at a
point $x \in A^{-1}(U)$ is equal to the Laplacian of $f$ at $A(x)$.

	In particular, $f \circ A$ is harmonic on $A^{-1}(U)$ if $f$
is harmonic on $U$.  On the complex plane, if $\phi$ is a holomorphic
function on an open set $V$ which takes values in another open subset
$U$ of ${\bf C}$, and if $f$ is a harmonic function on $U$, then $f
\circ \phi$ is a harmonic function on $V$.  If $f$ is also
holomorphic, then $f \circ \phi$ is holomorphic on $V$.  If $f(z)$
is a holomorphic function on an open subset $U$ of ${\bf C}$, then
\begin{equation}
	\widetilde{f}(z) = \overline{f(\overline{z})}
\end{equation}
is a holomorphic function on $\overline{U} = \{\overline{z} : z \in
U\}$.

	Let $A$ be an orthogonal linear transformation on ${\bf R}^n$
with associated matrix $(a_{j,k})$, so that the $j$th coordinate of
$A(x)$ can be written as
\begin{equation}
	\sum_{k=1}^n a_{j, k} \, x_k.
\end{equation}
We can associate to $A$ a linear isomorphism $\widehat{A}$ on the
Clifford algebra $\mathcal{C}(n)$, defined by the conditions that
$\widehat{A}(r) = r$ when $r$ is a real number,
\begin{equation}
	\widehat{A}(e_k) = \sum_{j=1}^n a_{j, k} \, e_j,
\end{equation}
and of course $\widehat{A}$ acting on a product of $e_l$'s is equal to
the product of the corresponding $\widehat{A}(e_l)$'s.  One can check
that the assumption that $A$ be an orthogonal transformation provides
the appropriate compatibility conditions, namely, $\widehat{A}(e_l)^2
= -1$ for $l = 1, \ldots, n$ and $\widehat{A}(e_l) \, \widehat{A}(e_m)
= - \widehat{A}(e_m) \, \widehat(e_l)$ when $l \ne m$, which basically
means that these products have real part equal to $0$.

	Suppose that $U$ is a nonempty open subset of ${\bf R}^n$, and
that $f(x)$ is a $\mathcal{C}(n)$-valued continuously differentiable
function on $U$, and let $f_A$ denote the function on $A(U)$ defined
by
\begin{equation}
	f_A(x) = \widehat{A}(f(A^{-1}(x))).
\end{equation}
One can check that $\mathcal{D}_L$ or $\mathcal{D}_R$ applied to $f_A$
at a point $x \in A(U)$ is the same as $\widehat{A}$ applied to
$\mathcal{D}_L$ or $\mathcal{D}_R$ of $f$, respectively, evaluated at
$A^{-1}(x)$.  In particular, $f_A$ is left or right
Clifford-holomorphic on $A(U)$ if $f$ is on $U$.

	It is easy to see from the definition of
Clifford-holomorphicity that at each point in the domain of a
Clifford-holomorphic function, the derivative of the function in the
direction of one of the coordinates is determined by the derivatives
of the function in the directions of the other coordinates.  Using the
invariance under orthogonal transformations described in the previous
paragraph, one can check that the derivative of a Clifford-holomorphic
function in any direction is determined by the derivatives in the
orthogonal directions.  Of course this is very clear in the classical
case of holomorphic functions of a single complex variable.

	We can also derive as follows.  Let $v$ be a vector in ${\bf
R}^n$, $\widehat{v} = \sum_{j=1}^n v_j \, e_j$ the corresponding
element of $\mathcal{C}(n)$, and $f(x)$ a $\mathcal{C}(n)$-valued
continuously differentiable function on an open subset $U$ of ${\bf
R}^n$.  If $\mathcal{D}_L f(x) = 0$ for some $x \in U$, then of course
\begin{equation}
	\widehat{v} \sum_{l = 1}^n e_l \,
		\frac{\partial}{\partial x_l} f(x) = 0,
\end{equation}
which can be rearranged as an identity between the directional
derivative of $f$ at $x$ in the direction of $v$ and a combination of
terms of the form $v_j \, \partial/\partial x_l f(x) - v_l \,
\partial/\partial x_j f(x)$ times $e_j \, e_l$.  These terms involve
directional derivatives of $f$ at $x$ in directions orthogonal to $v$,
and an analogous computation works for $\mathcal{D}_R f(x) \,
\widehat{v}$ when $\mathcal{D}_R f(x) = 0$.

\subsection{Some loose ends}
\label{section on some loose ends}

	Let $U$ be a nonempty open subset of ${\bf R}^n$.  If
$f_1$, $f_2$ are two real-valued continuously differentiable
functions with restricted support in $U$, consider the bilinear
form $B(f_1, f_2)$ acting on $f_1$, $f_2$ defined by
\begin{equation}
	B(f_1, f_2) = \int_U \langle \nabla f_1(x), 
					\nabla f_2(x) \rangle \, dx,
\end{equation}
where $\nabla f(x)$ denotes the gradient of $f$ at $x$, the vector
whose components are the partial derivatives
\begin{equation}
	\frac{\partial}{\partial x_j} f(x).
\end{equation}
Of course $B(f_1, f_2)$ is symmetric in $f_1$, $f_2$, and if $f_1$
is also twice continuously-differentiable, we have that
\begin{equation}
	B(f_1, f_2) = - \int_U  \Delta f_1(x) \, f_2(x) \, dx.
\end{equation}
Also notice that $B(f, f) \ge 0$, and indeed $B(f, f) = 0$ if and only
if $f = 0$, since we are considering functions with restricted
support.

	One can extend this to consider functions which tend to $0$ at
the boundary of $U$ and which have enough regularity, and it can be of
interest to make some assumptions about $U$ as well.  In particular,
for a nice domain $U$, one can try to diagonalize this bilinear form,
and thus the Laplacian, with respect to the standard integral innner
product
\begin{equation}
	\int_U f_1(x) \, f_2(x) \, dx,
\end{equation}
and for functions which tend to $0$ at the boundary of $U$.  This
means looking for smooth functions $f$ on $U$ which are eigenfunctions
of the Laplacian, so that there is a real number $\lambda$ such that
\begin{equation}
	\Delta f(x) = \lambda \, f(x)
\end{equation}
for all $x \in U$, and which tend to $0$ at the boundary of $U$.

	Just as for harmonic functions, one can define what it means
for a continuous function on $U$ to be an eigenfunction for the
Laplacian in the weak sense, and show that this implies that the
function is smooth.  There are also basic theories in analysis
for showing that there are eigenfunctions of the Laplacian which
tend to $0$ at the boundary, and enough of them to diagonalize the
Laplacian among functions which tend to $0$ at the boundary.
Of course there are some details in formulating this, and indeed
some variations along these themes.

	If $h(x)$ is a real-valued harmonic function on a nonempty
open subset of ${\bf R}^n$, and if $\phi(x)$ is a real-valued
continuously differentiable function with restricted support in $U$,
then
\begin{equation}
	\int_U |\nabla (h + \phi)(x)|^2 \, dx
		= \int_U |\nabla h(x)|^2 \, dx
			+ \int_U |\nabla \phi(x)|^2 \, dx,
\end{equation}
since
\begin{equation}
	\int_U \langle \nabla h(x), \nabla \phi(x) \rangle \, dx
		= - \int_U \Delta h(x) \, \phi(x) \, dx = 0.
\end{equation}
In particular, notice that
\begin{equation}
	\int_U |\nabla (h + \phi)(x)|^2 \, dx
		\ge \int_U |\nabla h(x)|^2 \, dx.
\end{equation}
More precisely, this works under suitable conditions for any
function $\phi$ which vanishes at the boundary of $U$.  Conversely,
one can use this as a basis for looking for harmonic functions
on $U$ with given boundary values, by trying to minimize 
\begin{equation}
	\int_U |\nabla f(x)|^2 \, dx
\end{equation}
among functions with those boundary values.

	Another approach to finding harmonic functions on a nonempty
open subset $U$ of ${\bf R}^n$ with specified boundary values is to
consider the supremum of subharmonic functions on $U$ which are less
than or equal to the prescribed values at the boundary.  A third
method is to generate a lot of harmonic functions on $U$ in general,
and then choose ones with desired properties.  We have already seen
``fundamental solutions'' on ${\bf R}^n \backslash \{0\}$ which are
harmonic or holomorphic and hence harmonic, and one can use translates
of these to get harmonic functions on ${\bf R}^n \backslash \{y\}$ for
each $y$ in ${\bf R}^n$, and then linear combinations of these with $y
\in {\bf R}^n \backslash U$ to get harmonic functions on $U$.

\end{document}